\documentclass{amsart}
\usepackage{ amsthm, amscd, amsfonts, amssymb}
\usepackage{hyperref}
\usepackage{amssymb}

\usepackage{amsmath}
\usepackage{amsfonts}

\theoremstyle{plain}

\newtheorem {lem}{Lemma}
\newtheorem {pr}{Proposition}

\newtheorem {Rem}{Remark}
\newtheorem {defi}{Definition}
\newtheorem {thm}{Theorem}
\newtheorem {corol}{Corollary}
\newtheorem {example}{Example}

\newcommand{\N}{\mathbb{N}}
\newcommand{\PP}{\mathcal{P}}
\newcommand{\Z}{\mathbb{Z}}

\newcommand{\XX}{\mathcal{X}}
\newcommand{\YY}{\mathcal{Y}}

\newcommand{\NN}{\mathcal{N}}

\newcommand{\C}{\mathbb{C}}
\newcommand{\inclust}{\varsubsetneq}

\newcommand{\inclu}{\subset}

\newcommand{\impliq}{\Rightarrow}

\newcommand{\be}{\begin{enumerate}}
\newcommand{\ee}{\end{enumerate}}
\newcommand{\bi}{\begin{itemize}}
\newcommand{\ei}{\end{itemize}}

\begin{document}

\title{On  Indecomposable  triples associated with nilpotent operators}

\author{A. Elkhantach}
\address{ A. Elkhantach, Center of Mathematical research in Rabat. Mohamed V  University in Rabat. Faculty of Sciences. BP 1014 Rabat Morocco}
\email{ }
\author{E. H. Zerouali}
\address{E. H. Zerouali,  Center of Mathematical research in Rabat. Mohamed V University in Rabat. Faculty of Sciences. BP 1014 Rabat Morocco}
\email{zerouali@fsr.ac.ma}
\thanks{ We are dept full to Professor 	Markus Schmidmeier  for his  enlightens some points in \cite{ring}.
 The Second author is  supported by the Project URAC 03 of the National center of research and by   Hassan II academy of sciences}

\begin{abstract}
We consider  in this paper the family of  triples $(V, T, U),$  where $ V$ is a
finite dimensional  space, $T $  is  a nilpotent linear operator on $V$ and  $U $ is an invariant  subspace of $T$. Denote $[U]= ker(T_{|U})$, and  $n_U= dim([U] )$.
 Our  main goal is to investigate possible   classification of  indecomposable  triples. The obtained classification
 depends on the order of nilpotency $p$,  on $n_U$ and on  $n_V$. Complete  classifications  are given for arbitrary $p$,  when  $n_U=1$, and when  $n_U=2$ and $n_V \le 3$. 
The case $ p \le 5$, treated in \cite{ring} is recaptured by using constructive proofs based on  linear algebra tools. The case $p\ge 6$, where the number of indecomposable triples is infinite, is also investigated.
\end{abstract}

\maketitle

\section{Introduction}
Let $V$  be a complexe vector space, $L(V)$ be the algebra of all linear operators on $V$  and  $T\in L(V)$.   A subspace $U$ of $V$ is said to be invariant for $T$ (or $T$-invariant) when $T(U) \subset U$. 
We will denote by $Lat(T)$ the lattice of all invariant subspaces for $T$.  Given any subset $A\subset V$, the smallest invariant subspace  generated by $A$ is 
$vect_T(A) = \{ \sum\limits_{k=1}^na_kT^kx, n \ge 1 , a_k\in {\mathbb C} \mbox{ and } x \in A\}$.
 For  $U\in Lat(T)$,
we  will say that  the triple $(V, T,U)$  is   indecomposable,
if for every two invariant subspaces  $V'$ and $ V''$, we have
$$\left\{\begin{array}{l}V = V' \oplus V''\\ U = (U \cap V') \oplus(U\cap V'') \end{array}\right.   \Longrightarrow 
 V' =\{0\}  \mbox{  or } \:  V'' = \{0\}.$$

We will also say that the invariant  subspace  $U$ is indecomposable in this case. \\

Recall that an operator is said to be algebraic if $P(T)=0$ for some polynomial $P$ and  is said to   be
  a nilpotent  of order $o(T)= p\in \N^*$ if $T^p=0$ and $T^{p-1}\ne 0$.\\
  Notice that  if $T$ is an algebraic  operator on an  infinite dimensional space, then there is  no indecomposable triple. This fact motivates our interset in the finite dimensional case. Also, since every algebraic operator is a direct summand of translations of nilpotent operators, we will restrict ourself to  nilpotent operators in finite dimensional spaces.\\

In the sequel $T$ is  nilpotent  of order $ p$. Obvious indecomposable triples are given by $({\mathbb C}^p,J_p, 0)$, where  $J_p$ is a cyclic Jordan block of order $p$. More generally, if $T$ has no reducing space, then every invariant subspace is indecomposable. On the other hand, normal operators, when $dim(V)>1$,  have no non trivial indecomposable subspaces.

 Our motivation in focusing on indecoposable triples is the next Krull-Remak-Schmidt property: 
Any triple is a direct sum of indecomposable triples and these direct summands
are unique up to isomorphism.

Richman-Walker in \cite{rich}  stated that if $T$ a nilpotent operator of order $5$,  admits    an  indecomposable triple  $(V, T,U) $, then necessarily $ dim( V)\leq 12$, $ dim(U)  \leq 6$, and $n_V \leq 3$
and these bounds are optimal. These observations, led to a complete characterization of indecomposable triples $(V,T,U)$   in  \cite{ring} in the case
 $ p  \leq 5$. See also \cite{bru} for more information  and details.\\

 Two triples  $(V, T,U)$ and $(V', T',U')$ are said to be similar if there exists  an invertible linear map
$f : V \longrightarrow V'$ with $f(U) = U'$ and such that  $fT = T'f$.  Similarity   induces   a classification of indecomposable triples that we will study below. More precisely, we will determine  similarity  classes of indecomposable  triples $(V, T,U)$ in the case $n_U=k$ with $k=1$ for arbitrary $n_V$  and  in the case $k= 2$ for  $n_V\leq3$   for  arbitrary $p$. We retrieve in particular several  results from\cite{ring} in the   case $p\le 5$. Our approach relies  on   direct computation of Jordan blocks and the notion of Jordan bases. It   does  not  require any notion from group theory.  The mains results are stated as follows :\\

\noindent {\bf Theorem A}: The number of the classes  of the  indecomposable triples  $(V,T,U)$ such that  $n_U=1$ is $C_{o(T)}^{2n_V-1}$.\\

\noindent {\bf Theorem B}:  The number of the classes  of the  indecomposable triples  $(V,T,U)$ such that  $n_U=n_V=2$,  is  $C_{o(T)}^4$.\\

\noindent {\bf Theorem C}: \begin{enumerate}
\item For $p\le 4$. There is no  indecomposable triple  $(V,T,U)$ such that  $n_U=2$ and $n_V=3$.
\item For $p=5$. There are exactly $8$  indecomposable triples  $(V,T,U)$ such that  $n_U=2$ and $n_V=3$.
\item For $p\le 5$. There are exactly $50$  indecomposable triples.
\item For $p\ge 6$. The number of the classes  of the  indecomposable triples  $(V,T,U)$ such that  $n_U=2$ and $n_V=3$ is infinite.
\end{enumerate}

\section{Basic tools}
\subsection{Height and valuation in  $v-$ modules.} 
By a \textbf{module} we will mean a module over a fixed discrete valuation
domain with prime $p$. A
\it{Valuated module} or \textbf{v-module},  is a module $B$ together
with a filtration $B=B(0)\supset B(1)\supset B(2)\supset\cdots$ such that
$p(B(n))\subset B(n+1)$.  If $x\in B(n) \setminus 
B(n+1)$ , we write $v(x)=n$  which is  called the \textbf{value} of $x$ is $n$. In the case where $B(n)=\{0\}$ for some minimal $n\ge 2$, the valuation will be said to be nilpotent of order $n$ and $B(n-1)$ will be called the kernel of the valuation. We will denote in this case $B(n-1)=[B]$.

Recall that a module $B$ is said to be a torsion module, if for every $x\in B$, there exists $n$ such that $p^nx=0$. The height $ht(x)$  is then defined as the minimum number $n$,  such that $p^nx=0$.  It follows  that $$[B]=\{x\in B \ : ht(x)=1 \}.$$

We state in the example below the context of  linear operators  which is our main area of investigations  in the sequel.
\begin{example} Let $V$ be a linear space and $T$ be a nilpotent operator on $V$ of order $p$. The space $V$ is   regarded as a  torsion  $v-$module over $\Z/p\Z$, under  the multiplication 
$$(\bar n, x) \in \\Z/p\Z\times V \to \bar n.x =: T^{n_0}x\in V. \ \ (n_0 \in \{ 0, \cdots, p-1\} \text{ and } \bar n= \bar n_0 ).$$
The filtration  $V(n)= T^n(V)$ induces a valuation on $V$ and  we  have 
$$V= V(0)\supset V(1) \supset \cdots \supset V(p-1)=ker(T) \supset  V(p) = 0.$$
  For every $x\in V$, we get $ht(x)= min\{n \mbox{ such that } T^nx =  0\}$   and 
$v(x) = max\{n : x \in T^n(V)\}$  $= max\{n \; ; \exists x_n \mbox{ such that } T^nx_n=x\}.$  In particular the valuation $v$ is nilpotent of order $p$.
 \end{example}

Notice in passing that an  invariant subspace $U,$  inherits  two natural   filtration from $V,$    $U(n) = T^n(U)$ and $U'(n)= U \cap V(n)$. Such observation is useful for possible classification on invariant subspaces.

 In contrast with the corresponding heights that are trivially  equal, the associated  value functions may be different. To deal with this fact, we adopt the  next definition from group theory, 
\begin{defi}
Let $x\in U$ be a non zero vector, the values of $x$ with respect to $U$, denoted here by
$ v_U(x),$   is defined as the unique number $ k \geq 0$ such that $x \in T^k(U) \setminus  T^{k+1}(U)$.  For convenience, we put  $v_U(0)=\infty$.
In the case where  $U=V,$   we will simply write  $v_V(x) := v(x)$ and is called the values of  x.
\end{defi}

We clearly have,  
\begin{itemize}
\item $v_U(x) \leq v(x)$  and $v_U(x) +ht(x) \le o(T_{|U}) $ for every $x \in U  \backslash \{0\}$.
\item $v_U\alpha (x+\beta y)\geq inf(v_U(x),v_U(y))$,  for every  $(x,y)\in U^2$  and  $\alpha, \beta$ nonzero numbers. In  the case where  $v_U(x)\neq v_U(y)$ the equality holds.
\item $v(\lambda x)=v(x)$ for $x \in U$ and $\lambda\in\mathbb{C}^*$.
\end{itemize}

We associate with $x\in U$, the notion of gap sequence  in $U$,  $gs_U(x)$ and the value  sequence  in $U$, $vs_U(x)$ that will play a crucial role in our approach. More precisely,
$$ \begin{array}{ll} gs_U(x) &=\{k\in [[0,q-1]] \ / v_U(T^{k+1}(x)) > v_U(T^k(x))+1\}\\ &
= \{k_1 <\cdots < k_r\},\\
vs_U(x) & =\{v_U(T^k(x)) \ /  \ k\in gs_U(x)\} \\
& = 
  \{v_1,\cdots,v_r\}.
\end{array}$$

The integer $r$ above will be called the length of  $x$.  We will drop the index $U$ in the sequel, in the case where $U=V$ or when there is no possible confusion.\\

    It is easy to see that   $k_r = ht(x) -1$,  $v_r = v_U(T^{k_r}x)$,  and that
  $v_{i+1} - v_i\ge 2$  for $i=1,\cdots,  r-1.$ Furthermore, we have the following useful properties on the gap sequence, 
\begin{pr}Let $T$ be nilpotent of order $p$, $U\in Lat(T)$ and $q= o(T_{|U})$. For every  $x\in U$, we have 
\begin{enumerate}
 \item $gs_U(Tx) \subset gs_U(x)-1 =  \{k_1-1 <\cdots < k_r-1\}$;
\item $ vs_U(Tx) \subset vs_U(x)$. Equality  holds  if and only if $k_1 \ne 0;$ 
\item $ card(gs_U(x))\le  min\{dim[V], k_r\};  $
\item $  2card(gs_U(x))-1 \le p.$ \end{enumerate}

\end{pr}
{\it Proof.} \begin{enumerate}
\item For every  $k\in gs_U(Tx)$  we have  $v_U(T^{k+1}(Tx))>v_U(T^{k}(Tx))+1$. It follows that  $v_U(T^{k+2}(x))>v_U(T^{k+1}(x))+1$ and hence that $k+1 \in gs_U(x)$. \item Is   clear.
\item  For the first assertion  $ card(gs_U(x))\le k_r$  is  immediate. For  $1\le i\le r-1$, we have $v( T^{k_i+1}x)=v_{i+1}-k_{i+1}+k_i+1$.We
consider  $b_{i+1}$ and $a_i$ such that $T^{v_{i+1}-k_{i+1}+k_i+1}b_{i+1}=T^{k_i+1}x,$  $T^{v_i}a_i= T^{v_{i+1}-k_{i+1}+k_i}b_{i+1}-T^{k_i}x$  and $a_r=b_r$.   From the expression $v_{i+1}=v(T^{k_{i+1}}x)=v(T^{v_{i+1}}b_{i+1})$, we deduce that $v_i=v(T^{k_i}x) < v_{i+1}-k_{i+1}+k_i=v( T^{v_{i+1}-k_{i+1}+k_i}b_{i+1})$, and hence that $v(T^{v_i}a_i)=v_i$. We have 
$$\begin{array}{lll} &i)&  \{T^{v_1}a_1,\dots,T^{v_r}a_r\}\inclu [V],\\ 
&ii)&  v(T^{v_i}a_i)=v_i \ \  for  \ \ 1\le i \le r,\\
 &iii) & v_1<\dots<v_r.
\end{array}$$ 

It follows that  the family $\{T^{v_1}a_1,\dots,T^{v_r}a_r\}$ is in $[V]$ and is  clearly independent. Thus  $ card(gs_U(x))\le  dim[V].$

\item Since $k_r = ht(x) -1$,    and $v_{i+1} - v_i\ge 2$  for $i=1,\cdots,  r-1.$  It follows that 
$$2(r-1)\le \sum\limits_{1}^{r-1}v_{i+1}-v_i \le v_r \le q-1.$$
\end{enumerate}
We introduce the next classical definition
\begin{defi} Let $V$ be a vector space and $T$ be a linear operator on $V$. We will say that $x\in V$ is $T-$consistent if $v(T^kx)=v(x)+k$ for every $k\le ht(x)-1$.
\end{defi} 
We clearly have: 
$$ \begin{array}{lll}
 x  \text{  is }   T-\text{consistent}  & \iff & gs(x)=\{ht(x)-1\}  \\ & \iff & vs(x) =\{v(x)+ht(x)-1\}.
 \end{array}
 $$
  We also have the following
\begin{lem}\label{cons}If $P=\sum\limits_{k=k_0}a_kX^k$ is a polynomial with $a_{k_0} \ne 0$ and $x\in V$. Then
\be \item  $v(P(T)x)=v(T^kx).$
\item If $k_0=0$, then $gs(x)=gs(P(T)x)$, and  $Vs(x)=vs(P(T)x)$.
\ee
In particular, if $k_0=0$, then  $x$ is $T-$ consistent if and only if $P(T)x$ is $T-$ consistent.
 \end{lem}
\subsection{ Values preserving property}
We  introduce  next the  notion of the  values preserving   property which is the main ingredient in our proofs.\\

\begin{defi}~~
 Let  $\YY := \{x_i\}_{i\in I}$ be  a family of vectors in $U$. We will say that $\YY$ satisfies the values preserving    property in $U$,  ( $VP_U$ for
 short and $VP$  in the case $ U=V$), if for every  finite subfamily $\{ x_j ;\: j\in J \subset I \} $, we have
$$ v_U(\sum_{i\in J  } \alpha_jx_j)= inf\{ v_U(x_{j})  : \mbox{ such that } \alpha_j \ne 0 \}.$$
\end{defi}
We have the next immediate observations
\begin{pr}\label{somme}\begin{enumerate} Let $V$ be a vector space and $T \in {\mathcal L}(V).$
\item If $V = U_1\oplus U_2$ with $U_1$ and $U_2$ are  invariant subspaces, then for every $x_1\in U_1$ and $x_2 \in U_2$, the famiy $\{x_1, x_2\}$ has $VP$.
\item  Let  $ \{x_i\}_{i\in I}$ be a family of vectors in $U$ and $\{P_i \: : i \in I\}$ be a family of polynomials such that $P_i(0)\ne 0$. Then, 

 $ \{x_i\}_{i\in I}$  satisfies  $VP_U$ if and only if  $\{P_i(T)x_i\}_{i\in I}$  satisfies  $VP_U$.
\end{enumerate}
\end{pr}
{\it Proof.}
\begin{enumerate}
\item Denote $v= v(\alpha x_1+\beta x_2)$ and consider $a \in V$ such that $T^v(a)= \alpha x_1+\beta x_2$. If we write $a=a_1 + a_2$ with $a_i \in U_i$, we will get  $ T^va_i\in U_i$ and hence $T^va_i= x_i$.  It follows that $v_{U_i}(x_i)\ge v$ for $i=1,2$
and then that $v= inf\{v_{U_1}(x_1), v_{U_2}(x_2)\}$.
\item  We write 
$$\begin{array}{ll}
 \sum\limits_{i\in J}\alpha_iP_i(T)x_i &= \sum\limits_{i\in J}\alpha_iP_i(0)x_i+\sum\limits_{i\in J}\alpha_i(P_i(T)-P_i(0))x\\ &  = \sum\limits_{i\in J}\alpha_iP_i(0)x_i+\sum\limits_{i\in J}\alpha_iQ_i(T)Tx_i, \end{array}$$ 
where  $Q_i(X) = \frac{P_i(X)-P_i(0)}{X}$ for $i \in I$.\\
 On the other hand, since $$v_U(\sum\limits_{i\in J}\alpha_iP(0)x_i) = inf\{v_U(x_i) : i\in J \}$$   and 
$$v_U( \sum\limits_{i\in J}\alpha_iQ_i(T)Tx_i)\ge  inf\{v_U(x_i)+1 : i\in J \},$$ 
we deduce that 
$$ \begin{array}{lll}
  inf\{v_U(P_i(T)x_i) : i\in J \} & =&   inf\{v_U(x_i) : i\in J \}\\
   &  =& v_U(\sum\limits_{i\in J}\alpha_iP_i(0)x_i)\\ 
   & =&v_U(\sum\limits_{i\in J}\alpha_iP_i(0)x_i+\sum\limits_{i\in J}\alpha_i(P_i(T)-P_i(0))x_i) \\ 
   &=& v_U(\sum\limits_{i\in J}\alpha_iP_i(T)x_i).
   \end{array}
$$
 Let     $\YY=\{y_1,\dots,y_r\}$ be a family  of non zero vectors  in $ U,$  and  denote   $n_i=ht(y_i)-1$. We introduce the next notations. 
\begin{itemize}
\item $[\YY]=\{T^{n_1}y_1,\dots,T^{n_r}(y_r)\}\subset [V], $ \item  $D(\YY)=\{T^{k}(y_i) \  \  /  \  \mbox 0\le i \le r \ {and} \ 0 \leq  k \leq n_i \},$
\item $ V(\YY) = span\{D(\YY)\}$. \end{itemize}
\end{enumerate}
The proof of the following lemma is easy and is left to the reader.
\begin{lem}\label{lem1} Let $\YY$ be a family of vectors in $ U.$  We have
\begin{enumerate}
\item If  $v_U(x) \ne v_U(y)$ for every  $x \ne y$ in $\YY$, then $\YY$ satisfies  $VP_U$.
 \item
If   $[\YY]$ satisfies $VP_U $, then $ D(\YY )$  is  linearly independent. In particular $D(\YY)$ is a basis  in $ V(\YY)$.
\end{enumerate}
 \end {lem}
We also have

\begin{lem}\label{lem2} Let $\YY$ be a family of vectors in $ U$ satisfying  $VP_U$. Then 
$$Tx \in V({\YY }) \iff x \in V({\YY }) +[V].$$
We deduce that, if $[V] \subset V({\YY })$, we have $$Tx \in V({\YY }) \iff x \in V({\YY }).$$ It follows in particular in this case that $V= V({\YY })$.
 \end {lem}

{\it Proof.} Only the direct implication requires a proof. Assume $Tx\in V(\YY)$ and write, 
$$Tx= a_1T^{n_1}x_1+ \cdots + a_kT^{n_k}x_k,  $$ with $x_i \in \YY$ and $a_i \ne 0$. From $VP_U$ if follows that $min\{n_i : i \le k\} \ge 1$ and then 
$$ x= a_1T^{n_1-1}x_1+ \cdots + T^{n_k-1}x_k + x -  (a_1T^{n_1-1}x_1+ \cdots + T^{n_k-1}x_k), \in   V({\YY }) +[V] .$$

 The last assertion follows from the fact that for every $x\in V$, there exists $q\ge1$ such that $T^q(x)\in[V]  \subset  V({\YY })$. The proof is complete.

  For $U \in Lat(T)$,   we  have 

\begin{pr}\label{lem3}
 Let ${\mathcal Y}$ be a family of vectors in $ U$ such that  $[\YY]$ is a basis in $[U]$. The following assertions are equivalent:
\begin{enumerate}
                                    \item $[\YY]$ satisfies $VP_U $;
                                    \item   $ D({\mathcal Y})$   is a basis  in  $U$.
                                  \end{enumerate}
\end{pr}

{\it Proof.}~ $(1)\Rightarrow (2)$: From  Lemma \ref{lem1},  we have $ {D(\mathcal Y})$  is an  independent family.

It remains to show   that  $ {D(\mathcal Y})$ spans  $U.$ To this aim,  let $x\in U$ be such that   $l= ht(x)$. We have  $$T^{l-1}x\in [V]\cap U =  span( [\YY])  \subset V({\mathcal Y}).$$
Using Lemma \ref{lem2},  we derive that $T^{l-2}x \in V({\mathcal Y }),$  and by  induction that $x\in V({\mathcal Y }).$

$(2)\Rightarrow (1)$.
Suppose  that   $ D({\mathcal Y})$   is a basis  in  $U$ and take $x\in  [U]$. Let  $\alpha_1,\ldots,\alpha_r$  be non-zero scalars and  $x_{i_1},\ldots,x_{i_r}\in [\YY]$ be such that $x=\alpha_1x_{i_1}+\cdots+\alpha_rx_{i_r}$
with $ v_U(x_{i_1})\leq...\leq  v_U(x_{i_r})$.
We claim that  $k=:v_U(x) \leq v_U(x_{i_1})$. Indeed, there exists non zero scalars $\beta_1,\ldots,\beta_r$  and $k_1,\ldots,k_r$ in $\mathbb{N}$ such that
$$
x=T^k(\beta_{i_1}T^{k_1}a_{i_1}+\cdots+\beta_{i_r}T^{k_r}a_{i_r}).
$$
Since $ {D(\mathcal Y})$  is a basis in  $U,$ we deduce that  $k+k_1=v_U(x_{i_1})$ and hence $k\leq v_U(x_{i_1})$.\\

We also show the next theorem to be used in  the sequel.

\begin{thm}\label{lem4}  Let  $\XX=\{a_1,\dots,a_r\}$  be  a family  of $T-$consistent vectors such that   $\{T^{hta_1-1}a_1,\dots,T^{ht(a_r)-1}(a_r)\}$ satisfies   $VP$.   Then

\begin{enumerate}
                           \item   $v(T^ka_i)=k$ for  $1 \le i \le r$ and $0\le k \le ht(a_i)-1$;
										\item $gs(a_i)= \{ht(a_i)-1\}$;
                           \item $D(\XX)$ has $VP$.

                         \end{enumerate}
                         \end{thm}
{\it Proof.} The first and the second assumptions are  trivial, we only show the last  one. To this goal,  consider
 $x=\alpha_1T^{k_1}a_{i_1}+\cdots+\alpha_sT^{k_s}a_{i_s}$ with $k_1\le k_2\le \cdots \le k_s$ and 
assume that $
v(\alpha_1T^{k_1}a_{i_1}+\cdots+\alpha_sT^{k_s}a_{i_s}) >k_1,
$
for some  $\{\alpha_1,\ldots,\alpha_s\}\inclu\mathbb{C}^*$.

 It will follow that there exists $P_1,\cdots,P_s$  polynomials in $\mathbb{C}[X]$  and  $l\in \N$ such that $k_1<l$ and
$$
T^l(P_1(T)a_{i_1}+\cdots+P_s(T)a_{i_s})=\alpha_1T^{k_1}a_{i_1}+\cdots+\alpha_sT^{k_s}a_{i_s}.
$$
Since $(T^ka_i)_{(i\in I,0\leq k \leq j_i)}$ are linearly independent and $l> k_1$,  we derive that $\alpha_1=0$.  Which is impossible.
\section{ Bases with values preserving  property. }

We devote this section to provide  an algorithm aiming to extend  independent families with values preserving property to  bases having     values preserving  property. The proofs of our main results  in  the remaining sections will rely heavily on such construction.\\

In what follows, 
$ U$ is  an invariant subspace of $T$ and $q = o(T_{\mid U})$.
For $n<m$ in $\N$, we denote  by $[[n,m]]=\{n,n+1, \cdots, m\}$ and for  $(l,k)\in  [[0,q-1]]\times[[0,p-1]]=:I(q,p)$,  we adopt the next notations

\begin{itemize}
\item  $U_{l,k}= [U]\cap T^l(U)\cap T^k(V)= \{x \in [U] \: \text{ : } v_U(x) \ge  l \text{ and } v(x)\ge k \}$, 
\item $U^{l,k}= U_{l,k}\ominus (U_{l+1,k}+U_{l,k+1})$,
\item  $\{x_i^{l,k}\}_{i\in I_{l,k}}$ a basis in $U^{l,k}$.
\end{itemize} 
It is clear that $[U]= U_{0,0}= \sum\limits_{(l,k) \in I(q,p)}U^{l,k}$.\\
 
We introduce now the next definition of mixed values preserving property,
\begin{defi}~~
Let $\mathcal{A}$ be a family of vectors in $[V]$.  We will say that $\mathcal{A}$ satisfies  $VP_{V,U}$  if it satisfies $VP$ and  $ \mathcal{A}\cap [U]$  satisfies  $VP_U$.
\end{defi}
We have the following
\begin{pr}\label{lem54}  Under the notations above, set $I = \cup _{(l,k)\in I(q,p)}I_{l,k}$ and ${\mathcal F}= \{x_i\}_{i\in I}$. Then 
 ${\mathcal F}$ is a basis in $[U]$ with  $VP_{V,U}$.
\end{pr}
\it{ Proof.}
 For $\{i_1,\dots,i_r\}\inclu I,$  $\{\alpha_1,\dots,\alpha_r \}\inclu \mathbb{C}^*$, and $x= \alpha_1x_{i_1}+\dots \alpha_rx_{i_r}$, we consider $l=\min \{v_U(x_{i_1}),\dots,v_U(x_{i_r})\}$and  $k=\min \{v(x_{i_1}),\dots,v(x_{i_r})\}$. 
 We  argue by induction on $n=p+q- (l +k+2)$.  \\
 $n=0$.   Since $(l,k) \in I(q,p)$,  we deduce from the expression $n=p-1-l +q-1-k=0$ that  $l=q-1$ , $k=p-1$  and hence  $v(x)=p-1$ and $v_U(x)=q-1$.\\

 Suppose now the result holds for $m=p +q -(l+k+2)-1$, and consider $n=p+q-(l+k+2)$.
 Let also $D =\{s\in[[1,r]] \ \ / \ \ v_U(x_{i_s})=l,\mbox{and} \ v(x_{i_s})=k\}$. We have the next two cases 
\begin{itemize}
\item  $D\neq\emptyset$. We will have  $x\notin U_{l+1,k}+ U_{l,k+1}$ and   then  $v(x)=k$, and
 $v_U(x)=l$.
\item $D=\emptyset$.   Let $D_1=\{s\in[[1,r]] \ \ / \ \ v_U(x_{i_s})=l\}$ and $D_2 =[[1,s]] \backslash D_1$. Denote $y=\Sigma_{s\in D_1}\alpha_sx_{i_s}$,
 and $z=\Sigma_{s\in D_2}\alpha_sx_{i_s}$.  We have $x=y+z$, and by induction hypothesis $v_U(y)=l$, $v(y)>k$, $v_U(z)>l$ and $v(z)=k$. Finally  $v(x)=k$, and
 $v_U(x)=l$. \end{itemize} The proof is complete.
\begin{flushright}
 $\blacksquare$
\end{flushright}

We derive the following  important result.
\begin{pr}\label{lem6}~
Let $V$ be a vector space and $T \in {\mathcal L}(V)$. We have 
 \begin{enumerate}
       \item Every  $VP$ family  in $[V]$  extends   to a basis in $[V]$ with $VP$.
       \item For every invariant subspace  $U$ there exists a basis $\mathcal{A}$ in $V$ verifying $VP_{V,U}$.
     \end{enumerate}

\end{pr}
\it{Proof.}~
It is obvious that  $(2)$  derives from $(1)$  and  Proposition \ref{lem54}. To prove $(1)$, let $\{x_1,\ldots,x_r\}$  be a  $VP$ family  in $[V]$  with $r < \dim [V]$, and $U_r=span\{x_1,\dots,x_r\}$.\\
It suffices to show that  there exists $x\in [V]$ such that $\{x_1,\ldots,x_r,x\}$  has $VP$.\\
 Since $r<dim([V])$,   there exists  $k\in \N$ such that
 $$[T^{k+1}(V)]\inclust U_r \inclust [T^k(V)].$$
For any given  $x\in [T^k (V)]\setminus U_r$, we have $v(x)=k$. To see that $\{x_1,\ldots,x_r,x\}$  has $VP$, let $y=\alpha_1x_{i_1}+\dots+\alpha_sx_{i_s}$ and $\alpha \in \C^*$. We have $k\le min\{v(x_{i_1}),\dots,v(x_{i_s})\}=v(y)$,  and hence, it suffices to show  that $v(y+\alpha x)= min (v(y),k)=k$. If $k< v(y+\alpha x)$, then $y+\alpha x \in  T^{k+1}[V]$. But $T^{k+1}[V]\inclu  U_r$, and $y\in U_r$, so $x\in U_r$ which is is not true.\\
For $W\in Lat(T)$, we denote  by $vs(W)=\{v(x)   /   x\in [W]  \ \mbox{and} \ x\neq 0\}$.
Let $\{x_1,\dots,x_n\}$  be  a basis in $[W]$ with $VP_W$. It is not difficult to see that  $$vs(W)=\{v(x_1),\dots,v(x_n)\}.$$ 

\section{Indecomposable triples  $(V,T,U)$ with $dim([U])=1$} We are concerned in this section with invariant subspaces $U$ satisfying $dim([U])=1$. Examples of such invariant subspaces are provided by cyclic invariant subspaces defined as follows: 
For $x\in V,$  the cyclic invariant  subspace generated by $x,$  given by $U_x =:span\{T^ix :\: i \ge
0\}.$ It is clear that $dim([U_x])=1$ and that every invariant subspace  $U$ such that $dim([U])=1$ is of the previous form.

We start with next structure  theorem, 
\begin{thm} \label{lem7} Let  $x\in V$ be a nonzero vector such that 
 $ gs(x) 
= \{k_1 <\cdots < k_r\}$ 
and  $
vs(x)=\{v_1 < \cdots < v_r\}.$
Then, there exists  $\{a_1,\dots,a_r \}\inclu V$ a family of $T-$ consistent vectors  such that
\begin{enumerate}
                                    \item  $ \{T^{v_i}a_i \ \ / \ \ 1 \leq i \leq r\}\inclu [V]$;                                    \item $x = T^{v_1-k_1}a_1+\dots+T^{v_r-k_r}a_r.$
                                  \end{enumerate}
\end{thm}
\it{Proof.}~~
We proceed by induction on the length $r$ of $x$.
For  $r=1$, we have  $gs(x)=\{k_1\}$, $vs(x)=\{v_1\}$,   $k_1=ht(x)-1$ and $v_1=v(T^{k_1}x)= v(x)+k_1$.  It follows that  $v(x)=v_1-k_1$.    If  $a_1\in V$  is such that $T^{v_1-k_1}a_1=x$, we  obtain  $T^{v_1}a_1\in [V]$ and $v(T^{v_1}a_1)=v_1$. The assertions is then proved for $r=1$.

\noindent Suppose  our assertion is true for $r-1$. Let $x$ be with  length $r$ and write  $gs(x)=\{k_1<\dots<k_r\}$.  Let now $\bar{x}=T^{k_1+1}x$;  we  clearly have $gs(\bar{x})=\{k_2-k_1-1,\dots,k_r-k_1-1\}$, and $vs(\bar{x})=\{v_2,\dots,v_r\}$.\\
 By our  induction hypothesis, there exists $\{a_2,\dots,a_r\} \inclu V$ such that 
\begin{enumerate}
                     \item $v(T^{v_i}a_i)= v_i; $ for $2\le i \le r$;
                                    \item  $ \{T^{v_i}a_i \ \ / \ \ 2 \leq i \leq r\}\inclu [V]$;
\item $T^{k_1+1}(x)=\bar{x} = T^{v_2-k_2+k_1+1}a_2+\dots+T^{v_r-k_r+k_1+1}a_r.$
                                  \end{enumerate}
It follows that $ x = T^{v_2-k_2}a_2+\dots+T^{v_r-k_r}a_r+y$, where $y\in Ker(T^{k_1})$. Since $$  v(x)  =v_1-k_1< v_2-k_2 \mbox{  and  }  v( T^{v_2-k_2}a_2+\dots+T^{v_r-k_r}a_r)  =v_2-k_2,$$  there exists  $a_1\in V$ such that $T^{v_1-k _1}a_1=y$. It is clear that $\{a_1,\dots,a_r\}$ is a family of $T-$ consistent vectors satisfying  $(1)$ and $(2)$.
\begin{Rem}~ Under the same notations as  before, we  write $ x_i=T^{v_i}a_i$. The family 
$\{x_1,\ldots,x_r\} \subset [V] $  has   $VP$  in $[V]$. If $r< \dim [V]$, it follows from Proposition \ref{lem54}, that there  exists  $\{x_{r+1},\dots, x_n\}\inclu [V]$ such that $\{x_1,\dots,x_n\}$  is a basis in $[V]$ with $VP$. We write $v(x_i)= v_i$ and we consider   $\{a_{1},\dots,a_n\}\inclu V$  such that  $x_i=T^{v_i}a_i$, for $1\le i \le n$. 
For  $\{b_1,\dots,b_r\}\inclu V$  a family of $T-$ consistent vectors  satisfying  $(1)$  and $(2)$ in Theorem \ref{lem4}, we set $z_i=T^{v_i}b_i$, for $1\le i \le r$. Let now  $\{z_{r+1},\dots, z_n\}\inclu [V]$ be  such that $\{z_1,\dots,z_n\}$ is a basis in $[V]$ with $VP$, and take  again $\{b_{r+1},\dots,b_n\}\inclu V$ such that  $z_i=T^{v_i}b_i$, for $r+1\le i \le n.$\\ It will  come  that  the families of vectors 
 $\{T^ka_i, 1\le i \le n \  \  \mbox{and} \  \   0\le k \le v_i\}$, and 
$\{T^kb_i, 1\le i \le n \  \  \mbox{and}  \  \ 0\le k \le v_i\}$ are a bases in $V$ with $VP$.\\
It is also easy to see that,
\begin{itemize}
\item $(V,T,vect_T\{a_1,\dots,a_r\})$ and  $(V,T,vect_T\{b_1,\dots,b_r\})$ are similar,
\item $(V,T,vect_T\{a_{r+1},\dots,a_n\})$ and  $(V,T,vect_T\{b_{1+1},\dots,b_n\})$ are  similar.
\end{itemize}

\end{Rem}

We denote by $
red_x=vect_T\{ a_1,\dots,a_r\}$ and $W(x)=vect_T\{ a_{r+1},\dots,a_n\}$. Then  $red_x$   and $W(x)$ are reducing subspaces such that  $V=red_x \oplus W(x)$.

\begin{thm} \label{lem8}Let $x, y \in V,$  and let $U_x=vect_T\{x\}$  and $U_y=vect_T\{y\}$ be the associated  cyclic invariant subspaces. Then
\begin{enumerate}
  \item $(V,T,U_x)$ is indecomposable if and only if $card(gs(x))=dim([V])$.
  \item The following are equivalent
\begin{itemize}
\item $(V,T,U_x)$ and $(V,T,U_y)$ are similar, \item $gs(x)=gs(y)$, \:$vs(x)=vs(y)$.\end{itemize}
\end{enumerate}

\end{thm}
{\it Proof.}~
\begin{enumerate}
  \item Suppose that  $(V,T,U_x)$ is indecomposable. Using the equality $V = red_x\oplus W(x)$,  and the inclusion $U_x\inclu red_x$, 
we derive that $W(x)=\{0\},$  and hence  $dim[V]=dim[red_x]$.
 Now, since $dim[red_x]=card( gs(x))$, we obtain    $$dim[V]= card(gs(x)).$$
      Conversely, suppose that  $card(gs(x))=dim([V])$ and let $V_1$ and $V_2$ in $ Lat(T)$, be such that $V=V_1 \oplus V_2$ and  $U_x=(V_1 \cap U_x)\oplus ( V_2 \cap U_x)$.  Since $dim[U_x]=1$,  either $V_1 \cap U_x=\{0\} $ or  $V_2 \cap U_x=\{0\}$. It follows that either  $U_x\subset V_1$ or  $U_x\subset V_2$. Suppose 
for  example that $U_x\inclu V_1$, then $V_2\cap U = \{0\}$. Using  $card(gs(x))=dim[red_x]=dim [V]$, and $red_x\inclu V_1$, we get  $[V]\inclu V_1$, and then  $V_2=\{0\}$, which implies that $(V,T,U_x)$ is indecomposable.
 \item Clearly,  if $(V,T,U_x)$ and $(V,T,U_y)$ are similar, then 
 $gs(x)=gs(y)$, and $vs(x)=vs(y)$.\\
For the other implication, suppose that 
 $gs(x)=gs(y)=\{k_1\d,\dots,k_r\}$, and that $vs(x)=vs(y)=\{v_1,\dots,v_r\}$. 
Since   $(V,T,U_x)$ is indecomposable, $r=card(gs(x))= dim[V]=n$.
There exists $\{b_1,\dots,b_n,c_1,\dots,c_n\}\inclu V$ such that:
  \begin{enumerate}
                     \item $v(T^{v_i}b_i)= v(T^{v_i}c_i)=v_i$ for $1\leq i \leq n$,
                     \item $x = T^{v_1-k_1}b_1+\dots+T^{v_n-k_n}b_n$ and $y = T^{v_1-k_1}c_1+\dots+T^{v_n-k_n}c_n$,
                      \item $\{ T^{v_1}b_1,\dots,T^{v_n}b_n\}$ and  $\{ T^{v_1}c_1,\dots,T^{v_n}c_n\}$ are bases in $[V]$ verifying $VP_V$.
                                 \end{enumerate}
 
Then, \\ $\{ T^kb_i \  \   /  \  \ 1\le i \le n   \  \  \mbox{and}   \  0\le k\le  v_i\}$  and 
 $\{ T^kc_i \  \   /  \  \ 1\le i \le n   \  \  \mbox{and}   \  0\le k\le v_i\}$ are basis in $V$ with $VP_V$. \\

The linear map defined on the basis by $\phi b_i=c_i$ is clearly an isomorphim from  $(V,T,U_x)$ to $(V,T,U_y)$.
\end{enumerate}
\begin{flushright}
 $\blacksquare$
\end{flushright}

Identifying similar cyclic indecomposable triples, we obtain the next result 
\begin{thm} \label{lem9} Let $T$ be a nilpotent operator of order $p$, then
The number of cyclic  indecomposable  triples  $(V,T,U)$  is $C_p^{2n_V-1}$.
\end{thm}
\it{Proof.}~~ Denote  $dim[V] =l$ and consider 
$$ {\mathcal P}=\{A\subset [[0,p-1]]  \ \  /  \  \  card(A)=2l-1\}.$$ 
We will show   that, the number of cyclic indecomposable triples $(V,T,U)$ with $n_V=l$   is exactly  $card({\mathcal P})$.

Let $(V,T,U)$ be an  indecomposable triple. There exists  a non zero vector  $x\in U$ such that $vect_T(x)=U.$
Since  $(V,T,U)$ is indecomposable, we get $card( gs(x))=card(vs(x))=l$. Let us  write $$vs(x)=\{v_1,\dots,v_l\}, \mbox{  and } gs(x)=\{k_1 <\dots< k_l\}.$$
It is not difficult to check  that $v_l= p-1,$ $  k_1\leq v_1$  and   for   $ 1\leq i\leq l-1$, we have 
$$ k_{i-1}<k_i  \mbox{  and  }  k_{i+1}-k_i < v_{i+1}-v_i.$$
It  follows that $$0\le k_1<  v_1+1 < v_2+1-(k_2-k_1)<v_2+1<\dots <v_l+1-(k_l-k_{l-1})<v_l+1=p.$$
In particular we will have $2l-1\leq p$.
We  denote by $$\phi(x)=\{k_1,v_1+1, v_2+1-(k_2-k_1),v_2+1 ,\dots ,v_l+1-(k_l-k_{l-1})\},$$  and we consider   the map $\phi$ that associates  with any indecomposable triple $(V,T,U_x)$, the set $\phi(x) \in {\mathcal P}$.

  From Theorem \ref{lem7}, we derive   that $\phi$ is well defined and is  one to one.
To see that $\phi $ is onto, we consider  $\{x_1<y_1<x_2<\dots<y_ {l-1}<x_l\}\in {\mathcal P}$, and we  denote  ${\mathcal V}=\{y_1,\dots,y_{l-1},y_l (=p-1)\}$.   By choosing an adequate Jordan form, we consider  $V$ be a vector space  such that  ${\mathcal V}= vs(V)$.

Let now $\{a_1,\dots,a_l\}\subset V$ be  such that,
$ht(a_i)=y_i$ for $1 \le i \le l$.   We set $x=\sum_{i=1}^lT^{y_i-k_i}a_i$, where   $k_1=x_1+1$, $k_i=k_{i-1}+y_i-x_i$ for $2\le  i \le l$.  We have $vs(x)=\{y_1,\dots,y_l\}$ and $gs(x)=\{k_1,\dots,k_l\}$. It follows that  $(V,T,U_x)$ is indecomposable, and  that $\phi(V,T,U_x)=\{x_1,y_1,x_2,\dots,y_{l-1},x_l\}$. Finally $\phi$ is onto and hence is bijective.
 We conclude that the number of the classes of the  indecomposable triples $(V,T,U)$ with $dim[V]=l$, $o(T)=p$, and $dim [U]=1$ is exactly $card({\mathcal P})=C_p^{2l-1}$.

  \begin{flushright}
 $\blacksquare$
\end{flushright}~

\section{Indecomposable triples  $(V, T,U)$ with $dim([U])= 2$} 

We study in this section   indecomposable triples $(V,T,U)$   such that $dim[U]=2$. Let  
 $\{x_1,\dots,x_n\}$ be a basis in $[V]$ with $VP_{V,U}$ and let  $\{x_{l_1},x_{l_2}\}= \{x_1,\dots,x_n\}\cap[U]$  be  a basis in $[U]$ with  $l_1 < l_2$. We write   $v(x_1) \le \cdots \leq v(x_{n-1})\leq  v(x_n)$ and we choose   $b_i\in V$ such that $T^{v(x_i)}b_i=x_i$. We have 
  $$ \{T^kb_i \ / \ \ 1\leq i \leq n, \mbox{and} \ 0\leq k \leq v(x_i)\}$$ is a basis in $V$ with $VP$.\\

  For $i \in \{1,2\}$, we denote  $n_i=v_U(x_{l_i})$, and  let $y_i\in U$ be  such that $T^{n_i}y_i=x_{l_i}$.  Similarly, we have 
    $$\{T^k y_i \ / \ \ 1\leq i \leq 2, \mbox{and} \ 0\leq k \leq n_{i}\}$$ is  a basis in $U$ with $VP_U$.

In the sequel, we write   $gs(y_1)=\{k_1<\dots<k_r\}$ and $vs(y_1)=\{v_1<\dots<v_r\}.$

\subsection{Construction of adequate  bases in  $V$} 

We start with some auxiliary observations to be used in our description.\\

By using Theorem \ref{lem7}, there exists  $\{a_1,\dots,a_r\}\inclu V$ such that:
\begin{enumerate}
                     \item $v(T^{v_i}a_i)= v_i; $ for $1\le i \le r$;
                                    \item  $ \{T^{v_i}a_i \ \ / \ \ 1 \leq i \leq r\}\inclu [V]$;
                                    \item $y_1 = T^{v_1-k_1}a_1+\dots+T^{v_r-k_r}a_r.$
                                  \end{enumerate}
The subspace  $red_{y_1}=vect_T\{a_1,\dots,a_r\}$ is a  reducing  subspace and since $v(x_{l_1})\le  v(x_{l_2})$, we get $x_{l_2}\notin red_{y_1}$. It follows that  there exists    $W(y_1)\in Lat(T)$ such that:
$$V= red_{y_1} \oplus W(y_1) \text{ and } x_{l_2}\in W(y_1).$$  For  $y_2=y_{2,1}+y_{2,2} \in red_{y_1} \oplus W(y_1),$  we denote $$gs(y_{2,2})=\{k_{r+1}<\dots<k_{r+s}\} \text{ and } vs(y_{2,2})=\{v_{r+1}<\dots<v_{r+s}\}.$$
 Applying  Theorem \ref{lem7} with  $y_{2,2}$, there exists  $\{a_{r+1},\dots,a_{r+s}\}\inclu W(y_1)$ such that:
\begin{enumerate}
                     \item $v(T^{v_i}a_i)= v_i; $ for $r+1\le i \le r+s$;
                                    \item $y_{2,2}= T^{v_{r+1}-k_{r+1}}a_{r+1}+\dots+T^{v_{r+s}-k_{r+s}}a_{r+s};$
                                    \item  $ \{T^{v_i}a_i \ \ / \ \ r+1 \leq i \leq r+s\}\inclu[W(y_1)]$.
  \end{enumerate}
Moreover, since $y_{2,1} \in red_{y_1} $ and $\{T^ka_i \:,  1\le i \le r \mbox \ ,  \  0\le k \le v_i\}$ is a basis in $red_{y_1}$, there exists $\{P_1,\dots,P_r\}\inclu \C[X]$ with $P_i(0)\ne 0$  and $\{l_1,\dots,l_r\}$ such that:\\
 \begin{equation}\label{eq1} y_{2,1}=P_1(T)(T^{l_1}a_1)+\dots+P_r(T)(T^{l_r}a_r).\end{equation}
 We have the following theorem 
\begin{thm}\label{lem10} Let $(V,T,U)$ be an indecomposable triple such that $n_U=2$. Then,
\be
 \item  $[T^{p-2}(V)]\inclu U$,
\item $v(T^{n_2}(y_2))=p-1$, and $\{p-1\}\inclust gs(y_2)$,
\item If moreover $dim[red_{y_1}]=n_V-1$, we can reduce to the case where $vs(y_2)=\{p-1\}\cup vs(y_{2,1})$.
\ee
\end{thm}
{\it Proof.}
\be
\item Let $y \in [T^{p-2}(V)]$ and write $y = \sum\alpha_iT^{s_i}a_i$. We have 
 $$p-2\le v(y)= v(\sum\alpha_iT^{s_i}a_i)=\inf\limits_{\alpha_i \ne 0} v(T^{s_i}a_i).$$
 Since $ v(T^{v_i}a_i))< v(T^{v_r}a_r)-1, $ for  $1\le i<r$ and 
	$ v(T^{v_{j}}a_{j}) <v(T^{v_{r+s}}a_{r+s})-1,$ for $r+1\le j < r+s$, we get 
 $$y= \alpha_rT^{v_r}a_r +  \alpha_{r+s}T^{v_{r+s}}a_{r+s}=   \alpha_rx_{l_1} +  \alpha_{r+s}x_{l_2}\in U.$$
\item From $(1)$, we have $v(T^{n_2}(y_2)=p-1$. Suppose that $vs(y_2)=\{p-1\}$, then $v(y_2)=p-1-n_2$.  If we choose  $a_n$ in such  way that $T^{p-1-n_2}(a_n)=y_2$, we get $(V,T,U)$ decomposable.
\item Since $dim[red_{y_1}]=n_V-1$, then $y_{2,2}= T^{p-1-n_2}a_n$.  If for some $1\le i \le r$, we have $p-1-n_2\le l_i$,   we will take  $a_n+ P_i(T)T^{l_i}a_i$   instead of $a_n$.
Hence for every $1\le i \le r$, if $P_i\neq 0$,  we get  $l_i< p-1-n_2$, and $v_i<p-1$.  Finally  $sv(y_2)=\{p-1\}\cup sv(y_{2,1})$.
\ee

The next theorem provides some sufficient  conditions for  a triple to be indecomposable.
\begin{thm}\label{lem11}
Under the notations  before, suppose that  there exists $\{d_1,d_2\}\inclu U$  and $\{k_1,k_2\}\inclu \N$  such that  :
\begin{enumerate}
                        \item $[red_{d_1}]+[red_{d_2}]=[V ]$,
                       \item $\{T^{k_1}y_1,T^{k_2}y_2\}$ has not $VP$.
                   \end{enumerate} Then,  $(V,T,U)$   is indecomposable if   one of the tree following conditions holds:
\begin{itemize}
\item[i)]  $0 < (k_1-k_2)((n_1-k_1)-(n_2-k_2))$,
                          \item[ii)] $v_1< v_2$, $\: \: \: k_2 < k_1$ and $ n_1-k_1=n_2-k_2$, 
                           \item[iii)] $k_1=n_1<k_2 < n_2$.
\end{itemize}

\end{thm}
 \it{Proof.}

It is obvious that  $iii)$ implies $i)$. We will show our assumption under the conditions $i)$ and $ii)$.  Seeking contradition, suppose that there are  non trivial subspaces $V_1$ and $V_2 $ in $Lat(T)$ such that 
$$V= V_1 \oplus V_2  \mbox{ and }  U = (V_1\cap U) \oplus (V_2\cap U).$$ We claim that 
$ (V_1\cap U) \ne \{0\} $ and $(V_2\cap U) \ne \{0\}$. Indeed, if for example  $U\subset V_1$, we will have 
 $[V]=[red_{d_1}]+[red_{d_2}]\inclu[V_1]$, and hence $V_2=\{0\}$  which is impossible.\\
It follows that,
 $$dim(V_1\cap [U ]) = dim(V_2\cap [U])=1.$$
  Consider  now $z_1\in V_1\cap U$ and $z_2\in V_2\cap U$ such that $vect_T\{z_1\}= V_1 \cap U$ and $vect_T\{z_2\}=V_2\cap U$. We deduce in particular that $U=vect_T\{z_1,z_2\}$.

Without loss of generality,  we assume that $ht(z_1)=ht(y_1)= n_1+1$ and $ht(z_2)=ht(y_2)=n_2+1$. We   write  
$$z_1 = P_1(T)y_1+ P_2(T)T^{l_1}y_2 \mbox{ and } z_2 = Q_1(T)T^{l_2}y_1+ Q_2(T)y_2,$$  where $P_i$,$ Q_i $ are polynomials  in $\C[X]$, such that  $P_i(0)\neq0$, $Q_i(0)\neq 0$ ,  $l_i\in\N$,  for $1\leq i \leq 2$, and  $n_2 - l_1 \le n_1\le l_2+n_2$.\\
Since  $\{T^{k_1}y_1,T^{k_2}y_2\}$ has not $VP$, we have $v(T^{k_1}y_1)=v(T^{k_2}y_2)( = v)$, and  there exists  $\alpha\in \C^*$ such that $v < v(T^{k_1}y_1+\alpha T^{k_2}y_2)$.\\
Suppose now that $i)$ is satisfied,  $0 < (k_1-k_2)((n_1-k_1)-(n_2-k_2))$. \\ Without loss of generality, we can assume that  $k_2<k_1$, and $ n_2-k_2<n_1-k_1$. It follows that $n_2<n_1$  and since $n_1\le n_2+l_2$, we deduce that $0<l_2$  and $k_1<n_1-n_2+k_2\le l_2+k_2.$
Moreover,  $$v(T^{k_1}z_1)=v(P_1(0)T^{k_1}y_1)=v=v(Q_2(0)T^{k_2}y_2)=v(T^{k_2}z_2).$$
On the other hand, we have $$ v< v(Q_2(0)P_1(0)[T^{k_1}y_1+\alpha T^{k_2}y_2])=v(Q_2(0)T^{k_1}z_1+\alpha P_1(0)T^{k_2}z_2).$$
From  $z_1\in V_1$,  $z_2\in V_2$, and $V=V_1\oplus V_2$,  by using Proposition \ref{somme}, we derive that $$v(Q_2(0)T^{k_1}z_1+\alpha P_1(0)T^{k_2}z_2)=min (v(T^{k_1}z_1),v(T^{k_2}z_2))=v. $$ Contradiction.\\
Suppose now that $ii)$ is satisfied,    $v_1< v_2$, $ k_2 < k_1$ and $ n_1-k_1=n_2-k_2$.\\  If  we assume $n_2+l_2=n_1$, it will  follow that  $T^{n_2}(z_2)=Q_1(0)T^{n_1}y_1+Q_2(0)T^{n_2}y_2$, and $T^{n_1}(z_1)=P_1(0)T^{n_1}y_1$. In particular $\{T^{n_1}z_1,T^{n_2}z_2\}$ will  not have $VP$. Which is again  a contradiction by using Propostion \ref{somme}.\\
Thus  $ n_1<n_2+l_2$, and hence $k_1=n_1+k_2-n_2<l_2+k_2$. We deduce as before  $$v(T^{k_1}(z_1))=v(P_1(0)T^{k_1}y_1)=v=v(Q_2(0)T^{k_2}y_2)=v(T^{k_2}z_2).$$ Since again we have, $$ v< v(Q_2(0)P_1(0)(T^{k_1}y_1+\alpha T^{k_2}y_2)=v(Q_2(0)T^{k_1}z_1+\alpha P_1(0)T^{k_2}z_2), $$ we obtain a contradiction.

 \begin{flushright}
 $\blacksquare$
\end{flushright}

   \subsection{Indecomposable triples $(V, T,U)$ when  $n_V\le 3$  and $n_U=2 $  } 
	For $y_1$ and $y_2$ given  as above,  there exist $\{a_1,\dots,a_n\}\inclu V$,  $\{s_1,\dots,s_l,r_1,\dots,s_k\}\inclu \N$, with $l+k\le 3$  and  $\{P_1,\dots,P_l\}\inclu \C[X]$ such that 
\be
\item $\{T^{v_1}a_1,\dots,T^{v_n}a_n\}$ is a basis in $[V]$ with $VP_{V,U}$,
\item $y_1=T^{s_1}a_1+\dots+T^{s_l}a_l$,
\item $y_{2.1}=P_1(T)T^{r_1}a_1+\dots+P_l(T)T^{r_l}a_l$,
\item $y_{2,2}=T^{r_{l+1}}a_{l+1}+\dots+T^{r_{l+k}}a_{l+k}$,
\item $P_1(0)\neq 0, \dots, P_l(0)\neq 0$.
\ee

We   simplify in a first step the expression of $y_{2.1}$.
\begin{pr}\label{lem12}Let  $(V,T,U)$  be an indecomposable triple such that $[U] =\{x_{l_1}, x_{l_2}\}$. Under the notations above,  without any  loss of generality, we can reduce to the following two  cases 
\begin{enumerate}
  \item  If  $card(vs(y_1))= 1$, or $card(vs(y_2))=2$, then 
$$
\left\{
\begin{array}{ll}  
y_1&=T^{s_1}a_1+\dots+T^{s_l}a_l,\\
y_{2,1}&=T^{r_1}a_1.
\end{array}
\right.
$$   
\item  If  $card(vs(y_1))= 2$, and $card(vs(y_2))=3$, then
  $$
\left\{
\begin{array}{ll}  y_1=P_1(T)T^{s_1}a_1+T^{s_2}a_2,\\
 y_{2,1}=T^{r_2}a_2+ T^{r_1}a_1.
\end{array}
\right.
$$       Where $P_1\in \C[X]$ is such that  $P_1(0)\neq0$.
  \end{enumerate}
\end{pr} 
 We need 
the next  auxiliary  lemma of independent interest.
\begin{lem}\label{lem13} Let  $a$ be  a non zero vector in $V$ and $Q\in \C[X]$ be such that $Q(0)\neq 0$. There exists $P\in \C[X]$ such that  $PQ(T)a=a.$
\end{lem}
{\it Proof.} Since $Q(0)\ne 0$, we have $ ht(b)=ht(a)$. Denote  $b=Q(T)a$, $ U_a=vect_T(a)$ and  $U_b=vect_T(b)$. We have $U_b\subset U_a$ and $dim(U_a)=dim(U_b)=ht(a)$. It follows that $U_a=U_b$ and hence there exists $P$ such that $a=P(T)b$.\\

{\it Proof.   $(1)$  If  $card(vs(y_1))= 1$, then $card(vs(y_{2,1}))= 1$. Also, if $card(vs(y_1))= 2$, then $n_V=3$, and hence  from Theorem \ref{lem10},  we deduce  that   $card(vs(y_{2,1}))= 1$.   We will discuss two sub-cases
  \begin{itemize}\item   $card(vs(y_1))= 1$. Because of cyclicity, we can write   $y_{2,1}=P_1(T)(T^{r_1}a_1)$,  with  $P_1\in \C[X]$, and $P_1(0)\neq 0$. From Lemma \ref{lem13}, there exists $P\in \C[X]$ such that 
$P(T)(y_{2,1})=P(T)P_1(T)T^{r_1}a_1= T^{r_1}a_1$.  We replace $y_2$ by $P(T)y_2$, to obtain  $y_{2,1}=T^{r_1}a_1$.\\
\item $card(vs(y_1))=2$. We have  $$vs(y_{2,1})\subsetneq vs(y_1)=\{v_1<v_2\}.$$
We consider $ j \in  \{1,2\}$ such that  $ht(y_{2,1})=v_j+1$ and $  i \in \{1,2\}\setminus \{j\}$.  Denote $r_j=v(y_{2,1})$, and let $\{b_j,b_i\}\inclu red_{y_1}$ be such that $$T^{r_j}b_j=y_{2,1},\  \mbox{and} \   \{T^{v_j}b_j,T^{v_i}b_i\} \ \ \mbox{is a basis in} \ \  [red_{y_1}].$$
 Since $v_1 < v_2,  the family $ $\XX=\{T^kb_l;  /  \  1\le l \le 2 \  \mbox{and} \ \ 0 \le k \le  v_l\}$ is a basis in $red_{y_1}$ with $VP$. It follows that $y_1$ can be written in the next form, $$y_1=Q_1(T)(T^{r_1}b_1)+Q_2(T)(T^{r_2}b_2),$$
where $\{Q_1,Q_2\}\inclu\C[X]$ are such that $Q_1(0)\neq 0$, and $Q_2(0)\neq 0$  and  with $\{r_1,r_2\}\inclu \N$.\\
By lemma \ref{lem13},  there exists $P\in\C[X]$ such that $$P(0)\neq 0,  \ \mbox{and} \  PQ_j(T)T^{r _j}b_j=T^{r_j}b_j.$$
Now, replacing $y_1$ by  $P(T)y_1$, and by setting $P_i=PQ_i$, we obtain   $$y_1=T^{r_j}b_j+ P_i(T)T^{r_i}b_i.$$
Also,  if we replace $a_j$  by $b_j$ and $a_i$ by $P_i(T)b_i$ respectively, we get 
$$y_1=T^{r_j}a_j+ T^{r_i}a_i \mbox{   and  }   y_{2,1}=T^{l_j}a_j.$$

\end{itemize}
(2) In the case where $card( gs(y_1))=2$, and $card(gs(y_2))=3$, we derive by using  Theorem  \ref{lem10} that  $card(gs(y_{2,1}))=2$. And because of $y_{2}\in red_{y_1}$,  $\{a_1,a_2\}$ can be chosen in such   way that  $$y_{2,1}= T^{r_1}a_1+T^{r_2}a_2,\  \mbox{and} \   \{T^{v_1}a_1,T^{v_2}a_2\} \ \ \mbox{is a basis in} \ \  [red_{y_1}], $$ It follows that $y_1$ is written as $$y_{1}= R_1(T)T^{s_1}a_1+R_2(T)T^{s_2}a_2,$$  where $\{R_1,R_2\}\inclu \C[X]$,  $R_1(0)\neq 0$, and $R_2(0)\neq 0$.\\
Again, there exists $P\in \C[X]$ such that 
$$P(0)\neq 0 \  \mbox{and} \  P(T)y_1=T^{s_2}a_2+PR_1(T)T^{s_1}a_1.$$ If we replace $y_1$ by $P(T)y_1$, and we take $P_1=PR_1$, we get  $$y_1=T^{s_2}a_2+P_1(T)T^{s_1}a_1.$$

As it has been shown  in Proposition \ref{lem12}, if $n_V=n_U=2$, then  there exists $\{a_1,a_2\}\inclu V$, $\{y_1,y_2\}\inclu U$, and $\{r_1,r_2,s\}\inclu \N$,  such that:
$$y_1=T^sa_1,\ \  and \ \ y_2=T^{r_2}a_2+ T^{r_1}a_1.$$

Since $\{n_2\} \inclust gs(y_2)$, we have $$gs(y_2)=\{v_1-r_1,n_2\}, v_1-r_1 < n_2 \  and \ r_1 < r_2.$$

 We derive  the next structure theorem, 
\begin{thm}\label{lem14}
Under the previous notations, let $(V,T,U)$ be a triple such that $ n_V=  n_U=2 $.  Then 
\begin{enumerate}
                       \item$(V,T,U)$ is indecomposable $\Leftrightarrow$  $n_1 < v_1-r_1;$
                       \item The triple $(r_1,r_2,s)$ is characterizes $(V,T,U)$;
                       \item There are $C_p^4$  indecomposable triples such that $ n_V=  n_U=2 $.  
                       \end{enumerate}
 \end{thm}
 \it{Proof.}
\begin{enumerate}
  \item It is sufficient to check that  the conditions in  the Theorem \ref {lem11} are satisfied.  It is clear that 
\begin{itemize}
 \item    $[red_{y_2}]=[V]$;
  \item  $\{T^{v_1-r_1}y_2,T^{n_1}y_1\}$ has not $VP$;
  \item $n_1<v_1-r_1<n_2$.
\end{itemize}
Conversly,  if $v_1-r_1 \leq n_1$, we will  get   $y= T^{r_2}a_2\in U$. From $y\in vect_T\{a_2\}$, $y_1\in vect_T\{a_1\}$, $U=vect_T\{y,y_1\}$ and $V=vect_T\{a_1,a_2\}$, we deduce  that  $(V,T,U)$ is decomposable. Which is a contradiction.
  \item We use  $s=v_1-n_1$, $r_2=v_2-n_2$, and $r_1=min\{v(x) \ \ / \ \ x\in U\}$.
\item  We have,  
$n_1 =v_1-s< v_1-r_1 < n_2 =v_2-r_2< v_2 - r_1$ and $v_2=p-1$. \\
Let $N_{p-1}=\{0,\dots,p-1\}$, $\NN_p=\{A\in \PP(N_p) \ \ / \ card(A)=4\}$, and $\AA$ the set of all class of triple $(V,T,U)$ indecomposable. The mapping $\phi$ defined by 
 $\begin{array}{cc}
  \phi : & \begin{array}{clll}
        \AA &\rightarrow&  \NN_p\\
        (V,T,U) & \mapsto & \{n_1,v_1-r_1,n_2,v_2-r_1\}, &
      \end{array}
\end{array}$  \\
 is bijective, and hence  $card(\AA)=card(\NN_p)=C_p^4$.

   \end{enumerate}
   \begin{flushright}
 $\blacksquare$
\end{flushright}
\subsection{Indecomposable triples $(V, T,U)$ with $n_V=3$ and $n_U=2$ }

 As before, we denote   $\{x_1,x_2,x_3\}$  for  a basis in $[V]$ with $VP_{V,U}$ property and 
for $1\le i \le 2$, we pose  $v_i=v(x_i)$,  $T^{v_i}a_i=x_i$  and $v_V=\{v_1, v_2, p-1\}$ with  $v_1\leq v_2\leq  p-1$.\\

In the case where  $(V,T,U)$ is   indecomposable,  by using  Theorem \ref{lem10}, we get $[T^{p-2}(V)]\inclu U$. In particular  $x_3\in U$, and $v_3 =: v(x_3)=p-1$.  On the other hand, since $ n_U=2$, either $x_1$ or $x_2$ belongs to $U.$  In the sequel, we  consider   $\{i,j\} = \{1,2\}$ such that $x_i\notin U$, and $x_j\in U$. We also write  $n_3=v_U(x_3)$, $n_j=v_U(x_j)$ and $y_2\in U$ such that 
$T^{n_j}y_2=x_j$. In particular,  we have $v_U=\{n_j,n_3\}$. Moreover,     if $x_1\in U$,  $vs(y_2) =  \{v_1\}$ and  if $x_2\in U$ , we obtain  $v_j \in vs(y_2) \subset \{v_1, v_2\}$.  

Let us  write as before $V = red_{y_2} \oplus W(y_2)$ with $W(y_2)\in Lat(T)$, and $x_3\in W(y_2)$. We also  have  $y_3=y_ {3,1} + y_{3,2}$  for some $y_ {3,1}\in red_{y_2}$ and  $y_ {3,2}\in Wy_2$.

\subsubsection{The case   $vs(y_2)=\{v_j\}$.}
We start with the next useful lemmas

\begin{lem}\label{lem15}
Let   $(V,T,U)$   be  indecomposable.  If $vs(y_2))=\{v_j\}$, then 
 $$max(r_1,r_2)< r_3 \mbox{ and } max(v_1-r_1,v_2-r_2)< v_3-r_3.$$
           \end{lem}
           \it {Proof.}
            Since $T^{v_3-r_3}y_3=x_3$, we get $max(v_1-r_1,v_2-r_2)< v_3-r_3$.  Let  now $k \in \{1, 2\}$, and suppose that $r_3\le r_k$. For $\tilde{a}_3= a_3+T^{r_k-r_3}a_k$, we have $T^{r_3}\tilde{a}_3=x_3$.  If $k=i$, $ y_3 =  T^{r_j}(a_j)+T^{r_3}\tilde{a}_3 $ and  if $k=j$,  $ y_3 =  T^{r_i}a_i+T^{r_3}\tilde{a}_3 $. In both cases $(V,T,U)$ will be decomposable. Contradiction.  
\begin{lem}\label{lem16}
Suppose that  $vs(y_2)=\{v_j\}$ and write  $   y_2=  T^{s_j}a_j\ \ \  \mbox{ and } \: \: \: \  \ \ 
         y_3 =  T^{r_1}a_1+T^{r_2}a_2+T^{r_3}a_3. $ Then 
         $(V,T,U)$   is decomposable if   one of the  following conditions holds:
         \be
         \item $r_i\le r_j$, and $v_j-r_j\le v_i-r_i$;
           \item $r_j \le r_i$, and $v_i-r_i < s_j-r_j$;
             \item $s_j\le r_j$.
           \ee
           \end{lem}
           \it {Proof.}
           \be
             \item Suppose that $r_i\le r_j$, and $v_j-r_j\le v_i-r_i$ and let us  consider  $\tilde{a}_i=a_i+T^{r_j-r_i}a_j$  and $\tilde{x_i}=T^{v_i}\tilde{a}_i$. If $v_j<v_i+r_j-r_i$, then $\tilde{x}_i=x_i$ and if $v_j=v_i+r_j-r_i$, then $i=1$, $j=2$, and $\tilde{x}_1=x_1+x_2$. In  both cases $\{\tilde{x}_i,x_j,x_3\}$ has $VP$. It is easy to check that  $y_3\in vect_T\{a_3,\tilde{a_i}\}$ and that  $y_2\in vect_T\{a_j\}$. Moreover
$$\left\{\begin{array}{l} V= vect_T\{a_3,a_i\} \oplus  vect_T\{a_j\},\\   U=U\cap  vect_T\{a_j\} \oplus U\cap vect_T\{a_3,\tilde{a_i}\}.\end{array}\right.$$ Thus $(V,T,U)$ is decomposable.
               \item   Suppose now  $r_j \le r_i$  and $v_i-r_i < s_j-r_j$. Let  us consider again $\tilde{a}_j=a_j+T^{r_i-r_j}a_i$.  We have  $$y_3= T^{r_3}a_3+  T^{r_j}\tilde{a}_j \ \mbox{ and} \  y_2=T^{s_j}\tilde{a}_j.$$
This yields $$U\inclu vect_T\{a_3,\tilde{a}_j\} \ \  and \ \ V= vect_T\{a_3,\tilde{a}_j\}\oplus vect_T\{a_i\},$$  and finally  $(V,T,U)$ is decomposable.\\

  \item  Consider $y=    T^{r_1}a_1+T^{r_3}a_3$. If  $s_j\le r_j$, then $y\in U$, and hence  $U=vect_T\{y,y_j\}$. Since moreover $y\in vect_T\{a_i,a_3\}$  and $y_2\in vect_T\{a_j\}$,  we get $(V,T,U)$ is decomposable.
               \ee
Denote in the  sequel $r_3=p-1-n_3$. We have 
\begin{thm}\label{lem17}Suppose that  $vs(y_2)=\{v_j\}$.  We have 
\begin{enumerate}
\item If  $(V,T,U)$ is indecomposable, then   
\begin{equation}\label{eqr}
 r_j <  s_j, \, max( r_i ,r_j)< r_3,  \:  \: and   \: max(v_i-r_i,   v_j-r_j)<p-1-r_3.
\end{equation}
Moreover, if we suppose that  \eqref{eqr} is satisfied, we obtain
\item If $vs(y_3)=\{ p-1>v_2 > v_1\}$, then 
 $(V,T,U)$ is indecomposable.

\item If  $vs(y_3)\neq\{ p-1>v_2 > v_1\}$  then,   $(V,T,U)$ is indecomposable, if and only if  $r_j\le r_i$, and $s_j-r_j\leq  v_i-r_i \le v_j-r_j$.

\end{enumerate}
 \end{thm}
\it{Proof.} 

\begin{enumerate}
\item  Since  $vs(y_2)=\{v_j\}$,  we have 
    $$   y_2=  T^{s_j}a_j \ \ \  \mbox{ and } \: \: \: \  \ \ 
         y_3 =  T^{r_1}a_1+T^{r_2}a_2+T^{r_3}a_3,  $$ 
where $s_j=v_j-n_j$ and $r_3=p-1-n_3$.
It is clear from Lemma \ref{lem15}, and Lemma \ref{lem16} that \eqref{eqr}  is satisfied.
 \item Suppose that $vs(y_3)=\{\ v_1 < v_2 <p-1\}$.   We get  $$[red_{y_2}]+[red_{y_3}]=[red_{y_3}]= [V] \ \mbox{and} \ \ v(T^{v_j-r_j}y_3)=v_j.$$ 
In addition, since  $\{T^{v_j-r_j}y_3, T^{n_j}y_2\}$ has not $VP$, we deduce that $n_j < v_j-r_j<n_3$. From Theorem \ref {lem11}, we derive   that $(V,T,U)$  is indecomposable.
\item Suppose first  that $r_i<r_j$. Since  $vs(y_3)\neq \{\ v_1 \le v_2 <p-1\}$, we get  $v_j-r_j\le v_i-r_i$.  Now,  Lemma \ref{lem16} says  that $(V,T,U)$ is decomposable,  and then   we have $r_j \le r_i$.

If $v_j-r_j< v_i-r_i$, then since $r_j\le r_i$, and   $vs(y_3)\neq \{\ v_1 , v_2 ,p-1\}$, we deduce that $r_i=r_j$. Now  by Lemma \ref{lem16}, the triple  $(V,T,U)$ is decomposable and then $v_i-r_i\le v_j-r_j$.

Suppose now  that $v_i -r_i<s_j-r_j.$  Again, by applying  Lemma \ref{lem16} $(V,T,U)$ is decomposable and hence $s_j-r_j\le v_i-r_i$.

Conversely, suppose that, $s_j-r_j \leq v_i-r_i \le v_j-r_j$, and $r_j \le r_i$.  It suffices to  show that  the conditions in Theorem \ref{lem11} are verified. Indeed, 
\begin{itemize}
\item Let $z=T^{s_j-r_j}y_3-y_2$.  We have $vs(y_3)=\{v_j,v_3\}$,  $vs(z)=\{v_i,v_3\}$,  $x_j\in red_{y_3}$, and $x_j \notin red_z$.  It follows that $[red_{y_3}]+[red_z]=[V]$.
\item $\{T^{v_j-r_j}y_3,T^{n_j}y_2\}$ has not $VP$  and  $n_j < v_j-r_j<n_3$.

\end{itemize}Finally  $(V,T,U)$ is indecomposable.

\end{enumerate}
\subsubsection{The case   $vs(y_2)=\{v_1,v_2\}$  and  $ vs(y_3)=\{v_1, p-1\}$ }  In the next theorem, the conditions   $v_3=p-1$ and $v_2\le v_3$  do not hold necessarily  in $(1)$.
\begin{thm}\label{lem18}
Assume that  $vs(y_3)=\{v_1,p-1\}$  and   $vs(y_2)=\{v_1,v_2\}$. We have
\be
\item If  $s_1\le r_1$ and $r_3\le s_2+r_1-s_1$, then   $(V,T,U)$ is decomposable if one of the next  two conditions holds
 \bi
\item[(i)] $n_2+s_1-r_1<n_3$,

\item [(ii)]  $n_2+s_1-r_1=n_3$, and $v_3\le v_2$.
   \ei
\item If $s_1< r_1$, then $(V,T,U)$ is indecomposable  if and only if $s_2+r_1-s_1- r_3<0$.
\item If $r_1< s_1$, then  $(V,T,U)$ is indecomposable if and only if  $n_3-s_1-n_2+r_1>0$.
\item  If $s_1= r_1$, then $(V,T,U)$ is indecomposable  if and only if $n_2<n_3$ and $s_2< r_3$.
\ee
\end{thm}
\it{Proof.}  
\be
\item Let $\tilde{y}_3=y_3-T^{r_1-s_1}y_2$, $\tilde{a}_3=a_3-T^{s_2-s_1+r_1-r_3}a_2$, and $\tilde{x}_3=T^{n_3}\tilde{y}_3$. Then,  if $n_2+s_1-r_1<n_3$, then $\tilde{x}_3=x_3 $ and  if $n_2+s_1-r_1=n_3$, we obtain $v_3\le v_2$, and $\tilde{x}_3=x_2+x_3$.\\
 In both cases, we get $\{x_1,x_2,\tilde{x}_3\}$ is a basis in $[V]$ with $VP_{V,U}$, and thus  $V=vect_T\{ \tilde{a}_3\} \oplus vect_T\{a_2,a_1\}$. Since in addition we have $\tilde{y}_3\in vect_T\{ \tilde{a}_3\}$, $y_2\in vect_T\{a_2,a_1\}$, and $U=vect_T\{\tilde{y}_3,y_2\}$, we deduce that $(V,T,U)$ is decomposable.
\item Suppose that  $s_2+r_1-s_1- r_3<0$. From $vs(y_3)=\{v_1,p-1\}$ and $vs(y_2)=\{v_1,v_2\}$, we get $[red_{y_3}]+[red_{y_2}]=[V].$\\ 
On the other hand, it is clair that  $\{T^{k_3}y_3,T^{k_2}y_2\}$, has not $VP$, with   $k_3=v_1-r_1$, and $k_2=v_1-s_1$.\\
Since $s_1< r_1$, $v_2\le v_3$, $n_2=v_2-s_2$ and $n_3=v_3-r_3$, we derive
        $$ \begin{array}{lll}
        0 &<&  (s_1-r_1)(s_2+r_1-s_1-r_3)\\ & \le &(s_1-r_1)(s_2-v_2+r_1-s_1+v_3-r_3)\\
          & = &(s_1-r_1)(n_3+r_1-n_2-s_1)\\ &=& (k_3-k_2)((n_3-k_3)-(n_2-k_2)).      
\end{array}$$
 By using Theorem \ref{lem11}, we get $(V,T,U)$ is indecomposable.\\

Conversely, suppose that    $(V,T,U)$ is indecomposable. Arguing by contradiction,  we assume  that $0\le s_2+r_1-s_1- r_3$. From  $v_2\le p-1$, we deduce
 $$0\le s_2-v_2+r_1-s_1+p-1- r_3=n_3-n_2+r_1-s_1.$$  If  $ n_3+r_1-n_2-s_1 = 0$, we obtain  $v_2=p-1$, and hence $(i)$ or $(ii)$ of Theorem \ref{lem18} $(1)$ is satisfied. It will follow that   $(V,T,U)$ is decomposable, which gives a contradiction.
\item  In the case  $n_3+r_1-n_2-s_1>0$, we  have $[red_{y_3}]+[red_{y_2}]=[V].$ Let $k_3=v_1-r_1$ and $k_2=v_1-s_1.$  It is clear that  $\{T^{k_3}y_3,T^{k_2}y_2\}$  has not $VP$.\\
Since moreover, $r_1< s_1$, we obtain  $$0<(s_1-r_1)(n_3+r_1-n_2-s_1)=(k_3-k_2)((n_3-k_3)-(n_2-k_2)), $$ 
and  by using Theorem \ref{lem11}, $(V,T,U)$ is indecomposable.\\
Conversely, suppose that    $(V,T,U)$ is indecomposable. Arguing by contradiction again,  suppose that $ n_3+r_1-n_2-s_1\le0.$ Since $n_3=p-1-r_3$ and $n_2=v_2-s_2$, we deduce that  $$0\le p-1-v_2\le r_3+s_1-r_1-s_2.$$  We have in addition  
$ n_3+r_1-n_2-s_1 \le 0$, and $v_2\le p-1$.  In particular,  one of the conditions in  $(1)$  of Theorem \ref{lem18}  is satisfied, and hence $(V,T,U)$ is decomposable. Contradiction.
\item Suppose that  $n_2<n_3$  and $s_2< r_3$. Seeking  contradiction, assume that $(V,T,U)$ is decomposable. Then  there exists $V_2$  and $V_3$ in $Lat(T)$ satisfying $$\left\{\begin{array}{ll} V& =V_2\oplus V_3, \\ U&=(V_2\cap U)\oplus (V_3 \cap U).\end{array}\right.$$ Since  $[red_{y_3}]+[red_{y_2}]=[V]$, we have  $$dim(V_2\cap [U ]) = dim(V_3\cap [U])=1.$$
  Consider  now $z_2\in V_2\cap U$ and $z_3\in V_3\cap U$ such that $vect_T\{z_2\}= V_2 \cap U$ and $vect_T\{z_3\}=V_3\cap U$. We deduce in particular that $U=vect_T\{z_2,z_3\}$.

Without loss of generality,  we assume that $ht(z_2)=ht(y_2)= n_2+1$ and $ht(z_3)=ht(y_3)=n_3+1$. We also  write  
$$z_2 = P_2(T)y_2+ P_3(T)T^{l_2}y_3 \mbox{ and } z_3 = Q_2(T)T^{l_3}y_2+ Q_3(T)y_3,$$  where  $l_i\in\N$ and  $P_i$,$ Q_i $ are polynomials  in $\C[X]$ such that $P_i(0)\neq0$, and $Q_i(0)\neq 0$.\\
Since $dim[V]=3$, we have $dim[V_2]=1$, or $dim|V_3]=1$, hence \\
$vs(z_2)=\{n_2\}$, or $vs(z_3)=\{n_3\}$.  It follows that $$v(z_2)=v_2-n_2=s_2 \ \ \mbox{or} \ \ v(z_3)=p-1-n_3=r_3.$$
On the other hand, we have $n_2<n_3$  and $ht(z_2)=n_2+1$, then $0<l_2$. Since  $v(y_2)=s_1=r_1=v(y_3)$, we obtain $v(z_2)=v(P_2(T)y_2)= v(y_2)=s_1<s_2$. Then  $v(z_3)=r_3$.\\
Moreover, since    $v(Q_3(T)y_3))=r_1<r_3$,  we obtain   $v(  Q_2(T)T^{l_3}y_2)=r_1$, and hence $l_3=0$.
We derive that $$z_3=Q_2(T)T^{s_2}a_2+ Q_2(T)T^{s_1}a_1+Q_3(T)T^{r_3}a_3+ Q_2(T)T^{r_1}a_1.$$ Thus $v(z_3)\le v(Q_3(T)T^{s_2}a_2)=s_2< r_3$. Contradiction.\\
Conversely suppose that $(V,T,U)$ is indecomposable.  Arguing by contradiction, assume  that $r_3\le s_2$, or $n_3\le n_2$. We use   condition $(i)$ or $(ii)$ in $(1)$  of Theorem \ref{lem18} $(1)$ to deduce that  $(V,T,U)$ is indecomposable.  Contradiction.
\ee

\subsubsection{The case   $vs(y_2)=\{v_1,v_2\}$ and $vs(y_3)=\{v_2,p-1\}$  } 
\begin{thm}\label{lem19}
Under the notations above, we suppose that  $vs(y_2)=\{v_1,v_2\}$ and $vs(y_3)=\{v_2,p-1\}$. Then, the  following properties hold
\begin{enumerate}
 \item If $s_2\le r_2$, then,  without any  loss of generality, we can reduce to the case 
     $  vs(y_2)=\{v_1,v_2\}$    and   $ vs(y_3)=\{v_1,p-1\}. $
 \item If $r_2< s_2$, then $(V,T,U)$ is indecomposable.
\end{enumerate}
\end{thm}
\it{Proof.}
\begin{enumerate}
\item For  $z_3=y_3-T^{r_2-s_2}y_2$, we have  $T^{n_3}(z_3)=x_3$, and $vs(z_3)=\{v_1,p-1\}$.
 \item  We have $[red_{y_3}]+[red_{y_2}]=[V]$,  $\{T^{v_2-r_2}y_3,T^{n_2}y_2\}$  not satisfying $VP$  and $n_2< v_2-r_2<n_3$. So  by using Theorem \ref{lem11},   we get $(V,T,U)$ is indecomposable.
\end{enumerate}

\subsubsection{ The case   $vs(y_2)=\{v_1,v_2\}$}~~ and $vs(y_3)=\{v_1,v_2,p-1\}$  

Let $\{r_1,r_2,r_3;s_1,s_2\}\inclu \N$, be such that $gs(y_2)=\{v_1-s_1<v_2-s_2\}$  and $gs(y_3)=\{v_1-r_1<v_2-r_2<v_3-r_3\}$. 
Recall that  $$\begin{array}{lll}
       y_2 &=&  P(T)  T^{s_1}a_1 +T^{s_2}a_2,\\
         y_3& = &T^{r_1}a_1+T^{r_2}a_2+T^{r_3}a_3.
      \end{array}$$
Let also $\{ \bar{a}_1, \bar{a}_2, \bar{a}_3\}\inclu V$, $\{ \bar{y}_2, \bar{y}_3\}\inclu U$, $\{ \bar{r}_1, \bar{r}_2, \bar{r}_3, \bar{s}_1, \bar{s}_2, \bar{s}_3\}\inclu \N$  and $ \bar{P}\in \C[X]$  be such that 
\be
\item $T^{v_2} \bar{a}_2=T^{n_2} \bar{y}_2$, and $T^{v_3} \bar{a}_3=T^{n_3} \bar{y}_3$.
\item $\{T^{v_1} \bar{a}_1,T^{v_2} \bar{a}_2,T^{v_3} \bar{a}_3\}$ is a basis in $[V]$ with $VP_{U,V}$.
\item  $gs(\bar{y}_2)=\{v_1-\bar{s}_1<v_2-\bar{s}_2\}$, and $gs(\bar{y}_3)=\{v_1-\bar{r}_1<v_2-\bar{r}_2<v_3-\bar{r}_3\}$.
\item $
       \bar{y}_2 =   \bar{P}(T)  T^{ \bar{s}_1} \bar{a}_1 +T^{ \bar{s}_2} \bar{a}_2,$ and $
          \bar{y}_3 = T^{ \bar{r}_1} \bar{a}_1+T^{ \bar{r}_2} \bar{a}_2+T^{ \bar{r}_3} \bar{a}_3.$
\ee
Since $\{T^{v_1} \bar{a}_1,T^{v_2} \bar{a}_2,T^{v_3} \bar{a}_3\}$ has $VP$, $T^{v_1}\bar{a}_1=\alpha_1x_1+x$, and  $T^{v_2}\bar{a}_2=\alpha_2x_2+z$, where $\{ \alpha_1,\alpha_2\}\inclu \C^*$, $x\in span\{x_2,x_3\}$ and $z\in span\{x_3\}.$\\

The main theorem in this case is stated as follows, 
\begin{thm}\label{lem20}
Under the previous notations, we have 
\begin{enumerate}
   \item The following are equivalent.\\
        $(i)$ For every  $ \{z_2,z_3\}\in U$, such that  $\{T^{n_1}z_2,T^{n_3}z_3\}$ is a basis in $[U]$  with $VP_{V,U}$, we have $v(z_2)=\{v_1<v_2\}$  and  $vs(z_3)=\{v_1<v_2<p-1\}$;\\
        $(ii)$ $r_1< s_1$,   $r_2< s_2$  and $s_1<r_1+ n_3-n_2$.
 \item If  $(i)$ is satisfied, then 
   \bi
\item (iii) $( \bar{r}_1, \bar{r}_2, \bar{r}_3, \bar{s}_1, \bar{s}_2,\frac{\alpha_1}{\alpha_2} \bar{P}(0))=(r_1,r_2,r_3;s_1,s_2,P(0)),$
\item (iv) $(V,T,U)$ is indecomposable.
\ei

\end{enumerate}
 \end{thm}
\it{Proof.} \be

\item  $(i)\Rightarrow (ii).$
Let $i\in \{1,2\}$. Suppose that $s_i\le r_i$, and   consider $y=y_3-T^{r_i-s_i}y_2$. We have $v(T^{n_3}y)=p-1$, and $v_i \notin vs(y)$.\\
If $n_3-n_2+ r_1\le s_1$, then if we take  $y=y_2-T^{s_1-r_1}y_3$. We have $\{T^{n_2}y,T^{n_3}y_3\}$ is a basis in $[U]$  with $ VP_{V,U}$, and $ vs(y)=\{v_2\}$.\\
 $(ii)\Rightarrow (i)$ Let $ \{z_2,z_3\}\in U$ be such that  $\{T^{n_2}z_2,T^{n_3}z_3\}$ is a basis in $[U]$  with $VP_{V,U}$. There exists $P_2$, and $P_3$ in $\C[X]$ such that $z_3=P_2(T)y_2+P_3(T)y_3$. Since $r_2 < s_2$, and  $r_1< s_1$,   we have $P_3(0) \neq 0$, and  
  $$\begin{array}{ll}v(T^kP_3(T)y_3)& =v(P_3(T)T^ky_3)=v(T^ky_3) \\ & < v(T^{k}y_2)\le v(P_2(T)T^{k}y_2)\\ &= v(T^kP_2(T)y_2),\end{array} $$  for every $k\in \N$.
  It follows that  $v(T^kP_3(T)y_3) <  v(P_2(T)T^{k}y_2)$. Thus $$ \mbox{for every }  k\in \N, \ v(T^kz_3)=v(T^kP_3(T)y_3)=v(T^ky_3),$$
and then   $$gs(z_3)=gs(y_3) \: \  and  \: vs(z_3)=vs(y_3).$$
Also, there exists $Q_2$ and $Q_3$ in $\C[X]$ such that $z_2= Q_2(T)y_2+Q_3(T)T^{n_3-n_2}y_3$.  Since $s_1<r_1+n_3-n_2$, $s_2<r_2+n_3-n_2$ and $\{T^{n_2}z_2,T^{n_3}z_3\}$ has $VP$,   we have $Q_2(0)\neq 0$  and \\
  $\begin{array}{ll}v(T^kQ_2(T)y_2)& =v(Q_2(T)T^ky_2)=v(T^ky_2)   < v(T^{k+n_3-n_2}y_3)\\ & \le v(Q_3(T)T^{k+n_3-n_2}y_3) = v(T^kQ_3(T)T^{n_3-n_2}y_3),\end{array} $\\
  for every $k \le n_2$.
  It follows that  $v(T^kQ_2(T)y_2) <  v(Q_3(T)T^{k}T^{n_3-n_2}y_3)$. Thus, for every  $  k\in \N,$ we have $$ v(T^kz_2)=v(T^kQ_2(T)y_2)=v(T^ky_2),$$
and then   $$gs(z_2)=gs(y_2) \: \  and  \: vs(z_2)=vs(y_2).$$ 
\item
\bi
\item $(iii)$ From the proof  of $(i)\impliq (ii)$, we obtain
 $$\begin{array}{l} \{v_1-\bar{s}_1,v_2-\bar{s}_2\}=gs(\bar{y}_2) =gs(y_2)= \{v_1-s_1,v_2-s_2\}\\ \{v_1-\bar{r}_1,v_2-\bar{r}_2,v_2-\bar{r}_2,v_3-\bar{r}_3\}=gs(\bar{y}_2)=gs(y_2) =\{v_1-r_1,v_2-r_2,v_3-r_3\}.\end{array}$$

 It follows that  $( \bar{r}_1, \bar{r}_2, \bar{r}_3, \bar{s}_1, \bar{s}_2)=(r_1,r_2,r_3,s_1,s_2)$.\\
Now, let us  show that $\frac{\alpha_1}{\alpha_2} \bar{P}(0)=P(0)$.\\
There exist $\{P_1,P_2,P_3,Q_1,Q_2,Q_3,R_1,R_2,R_3\}\inclu \C[X]$ such that 
 $$\begin{array}{lll}
         \bar{a}_1&=& P_1(T)a_1+P_2(T)T^{v_2-v_1}a_2+P_3(T)T^{v_3-v_1}a_3, \\
                           \bar{a}_2&=& Q_1(T)a_1+Q_2(T)a_2+Q_3(T)T^{v_3-v_1}a_3,\mbox{ and } \\
          \bar{a}_3&=& R_1(T)a_1+R_2(T)a_2+R_3(T)a_3.\\
     \end{array}$$
We have\\$\begin{array}{lll}
         \bar{y}_2&=& T^{s_2}\bar{a}_2+\bar{P}(T)T^{s_1}\bar{a}_1\\
                          &=& Q_1(T)T^{s_2}a_1+Q_2(T)T^{s_2}a_2+Q_3(T)T^{v_3-v_1+s_2}a_3\\
                         & &+ \bar{P}P_1(T)T^{s_1}a_1+\bar{P}P_2(T)T^{v_2-v_1+s_1}a_2+\bar{P}P_3(T)T^{v_3-v_1+s_1}a_3\\
                 &=& [Q_1(T)T^{s_2}+ \bar{P}P_1(T)T^{s_1}]a_1+[Q_2(T)T^{s_2}+\bar{P}P_2(T)T^{v_2-v_1+s_1}]a_2\\
                    & &+ [Q_3(T)T^{v_3-v_1+s_2}+\bar{P}P_3(T)T^{v_3-v_1+s_1}]a_3.\\\\
               \end{array}$\\
  We deduce that 
  
  $ T^{v_1-s_1} \bar{y}_2= Q_2(T)T^{v_1-s_1+s_2}a_2+\bar{P}(T)(P_1(0)x_1+P_2(0)x_2  +P_3(0)x_3) $
and 

$ Q_2(0)T^{v_1-s_1} y_2 = Q_2(0) P(0)x_1+Q_2(0)T^{v_1-s_1+s_2}a_2.$

In particular $$[P_1(0)\bar{P}(T)- Q_2P(0)]x_1+[Q_2(T)-Q_2(0)]T^{v_1-s_1+s_2}a_2\in U.$$ From $T^{v_1-s_1+s_2+1}a_2\in U$, we get  $[Q_2(T)-Q_2(0)]T^{v_1-s_1+s_2}a_2\in U$, and hence $P_1(0)\bar{P}(T)x_1- Q_2(0) P(0)x_1\in U$. But since  $x_1\notin U$,  we derive that $P_1(0)\bar{P}(T)-\bar{P}(T) Q_2(0) P(0)=0$.
It easy to see that $\alpha_1=P_1(0)$, and $\alpha_2=Q_2(0)$. Finally  $\bar{P}(0)\alpha_1=P(0)\alpha_2$.
\item $(iv)$We have $[V] = [red_{y_2}]+ [red_{y_3}]$, $v(T^{v_2-r_2}y_2)= v_2\in vs(y_2)$, and $r_j < s_j$. By using Theorem \ref{lem11} again, we obtain $(V,T,U)$ is indecomposable. 

\ei
\ee

Notice that if  condition $(1)$ of the previous theorem is fulfilled, then  
$$ r_1 < s_1 \le v_1<v_2-1 <v_2 <v_2+1 < v_3=p-1,$$
which forces $p\ge 6$. Conversely using the previous theorem, we  deduce the next structure corollary
\begin{corol} Let  $ T$  be nilpotent of order $p$.  Then the  number of indecomposable triple such that $n_V=3$ and $n_U=2$ is infinite if and only if $p\ge 6$.
\end{corol}

{\it Proof.} From  $(2)- ii) $ in Theorem \ref{lem20}, for isomorphic  indecomposable triples  we have  $$( \bar{r}_1, \bar{r}_2, \bar{r}_3, \bar{s}_1, \bar{s}_2,\frac{\alpha_1}{\alpha_2} \bar{P}(0))=(r_1,r_2,r_3;s_1,s_2,P(0)).$$ It follows  then that the indecomposable triples contains an infinite familly ${\mathbb C}$-indexed.

\section{Determination of  indecomposable triples $(V, T,U)$ with Nilpotency Index $\le  5$}
We  suppose  that $o(T) \leq 5 $. Our main objective in this section is to exhibit  all  non isomorphic indecomposable triples. It is known that, in this case  such triples exist if and only if  $n_V \leq 3$, see \cite{bru} for example. Moreover, we will see that  necessarily, we have $n_U<n_V$. For this reason, we s we will restrict ourself first,  to the case  $ n_U< n_V \leq 3$ .\\
 Let $\{x_1,x_2,x_3\}$ be a basis in $[V]$ with $VP_{V,U}$ and  let $v(x_i)=v_i$ for $1\leq i \leq 3$ be  ordered such that  $v_1 \leq v_2\leq p-1$.

\begin{thm}\label{lem21}
  Let $ p\le 5$ an $(V,T,U)$ be a triple, then 
\begin{itemize}
\item  If $ n_V=3$, and $n_U=2$ and  $(V,T,U)$  is indecomposable then  $p=5$.
\item 
The is  exactly  8 indecomposable triples  $(V,T,U)$  such that $p=5$, $ n_V=3$, and $n_U=2$.
\end{itemize}
\end{thm}
\it{Proof.}~ We notice first that if  $(V,T,U)$ is indecomposable, then by using Theorem \ref{lem10}, we get $x_3\in [U]$ and  $1<card(vs(y_3))$. It follows that  $1\le r_3\le p-2\le 3$, and $1\le n_3\le3$. We will determinate all indecomposable  triples $(V,T,U)$. We distinguish the next alternative cases,
\begin{itemize}
\item[$(a)$]   For every  $  y\in U,$  such that  $v(T^{n_3}y)=p-1$, we have $$vs(y)=\{v_1,v_2,p-1\};$$  \item[$(b)$] $vs(y_3)\neq \{v_1,v_2,p-1\}$, and $vs(y_2)=\{v_j\}$; 
\item[$(c)$] $vs(y_3)=\{v_1<p-1\}$, and $vs(y_2)=\{v_1<v_2\}$; \item[$(d)$] $vs(y_3)=\{v_2<p-1\}$, and $vs(y_2)=\{v_1<v_2\}$.\end{itemize}
$(a)$ Suppose that for every $  y\in U,$  such that  $v(T^{n_3}y)=p-1$, we have $vs(y)=\{v_1,v_2,p-1\}$. We will get in particular $vs(y_3)=\{v_1<v_2<p-1\}$. We write $gs(y_3)=\{k_1<k_2<k_3\}$, then clearly $$0\le k_1\le v_1 <v_2-k_2+k_1<v_2<p-1-k_3+k_2< p-1\le4.$$
It follows that $p=5$, $k_1=v_1=0$, $v_2=2$, $k_2=1$ and $k_3=2$.  Since $gs(y_3)=\{v_1-r_1,v_2-r_2,p-1-r_3\},$ we deduce that $r_1=0$, $ r_2=1$ and $r_3=2$.   If $vs(y_2)=\{v_1,v_2\}$, then  since $v_1=0$, $v(T^{p-1}(y_3-y_2))=p-1$, and $vs(y_3-y_2)\inclu \{v_2,v_3\}$. Thus $vs(y_2)=\{v_j\}$. Also from  $v_1=0$, we deduce $x_1\notin U,$  and $x_2\in U,$  in particular $j=2$.  Using  Theorem \ref{lem17}, we get $1=r_2< s_2\le v_2=2$.Thus $s_2=2$, and finally,  we derive that there in only one indecomposable triple in this case,   $$(1) \ \  \ y_2=Ta_2, \  y_3=T^2a_3+Ta_2+a_1, \  (v_1,v_2,v_3,n_2,n_3)=(0,2,4,0,2).$$
$(b)$ Suppose that  $vs(y_3)\neq \{v_1,v_2,  p-1\}$ and $vs(y_2)=\{v_j\}$. It will follow that   $vs(y_3)=\{v_j,v_3\}$.
By Theorem \ref{lem17},  we have $$r_j\le r_i<r_3, \ 0< s_j-r_j\le v_i-r_i \le v_j-r_j<p-1-r_3<p-1\le 4. $$  Thus   $2 \le 4-r_3$ and then  $r_3\in \{1,2\}$.
\bi\item $r_3=1.$ From  $r_1<r_3$ and $r_2< r_3$, we derive that $r_1=r_2=0$ and that $0<s_j \le v_i\le v_j\le2$. Using $r_i=r_j$ and   Lemma \ref{lem16},  we get $v_i-r_i<v_j-r_j$. We deduce  that $0<s_j \le v_i< v_j\le2$. So $i=1$, $j=2$,  $1=s_2=v_1$, and $v_2=2$. It follows that $2= v_j-r_j<p-1-r_3<p-1\le 4$, and in particular  $p=5$.  Hence
$$(2) \  \ \ y_2=T^2a_2, \  y_3=T(a_3)+a_2+a_1, \  (v_1,v_2,v_3,n_2,n_3)=(1,2,4,0,3).$$
\item $r_3=2.$  We have $r_1 \le 1$, $r_2\le1$, $ p-1-r_3\le 2$, $ v_1-r_1\le 1$, and $v_2-r_2\le 1$. Since $0<s_j-r_j\le v_i-r_i\le v_j-r_j$ and $r_j\le r_i$, we deduce that $v_j-r_j=v_i-r_i=s_j-r_j=1$ and  that $v_j\le v_i$. If $r_j=r_i$, we will get  $v_i=v_j$,  and $(V,T,U)$ will be decomposable. Hence $r_j=0$, $ r_i=v_j=1$, $ v_i=2$, $ j=1$ and $s_j=1$. From  $1=v_j-r_j<p-1-r_3= p-1-2\le 2$, we deduce that $p=5$. Finally we get $$(3) \ \  \ y_1=Ta_1, \  y_3=T^2a_3+Ta_2+a_1, \ and  (v_1,v_2,p-1,n_1,n_3)=(1,2,4,1,2)$$
\ei
$(c)$ Suppose that $vs(y_2)=\{v_1<v_2\}$ and $ vs(y_3)=\{v_1<p-1\}$. We distiguish three cases
\bi\item  $r_1=s_1$. We have  $s_2<r_3$  and $v_2-s_2<p-1-r_3$. Thus $v_2+1<p-1$,  and since  $vs(y_2)=\{v_1< v_2\}$, we derive that 
$$  0\le s_1< s_2< v_2 \mbox{ and } 0\le v_1<v_1+1<v_2<v_2+1<p-1=4.$$
Hence $v_1=r_1=s_1=0$, $ s_2=1$, $ v_2=2$, $ p=5$ and $r_3=2$.
 Finally we get $$(4) \ \ \ y_2=a_1+Ta_2, \ y_3=T^2(a_3)+a_1, \ (v_1,v_2,p-1,n_2,n_3)=(0,2,4,1,2).$$
\item $r_1< s_1$. We obtain by using  Theorem \ref{lem18},   $n_3+r_1-n_2-s_1>0$. Since moreover $vs(y_3)=\{v_1<p-1\}$ and $vs(y_2)=\{v_1<v_2\}$,  we get \\
$0 \le v_1-s_1< v_2-s_2=n_2< n_3-s_1+r_1=p-1-r_3-s_1+r_1<p-1-s_1<p-1\le 4$. Hence $p=5$, $p-1-s_1=3$, $n_3-s_1+r_1=2$, $ v_2-s_2=1$ and $ v_1-s_1=0$. It follows that  $v_1=s_1=1$,  $r_1=0$, $n_3=3$, $r_3=1$, $n_2=1$ and $(s_2,v_2)=(2,3)$ or $(s_2,v_2)=(3,4)$. Finally we get 
$$(5) \ \ \ y_2=Ta_1+T^2a_2, \ y_3=Ta_3+a_1, \  (v_1,v_2,p-1,n_2,n_3)=(1,3,4,1,3),$$or 
$$(6)  \ \ y_2=Ta_1+T^3a_2, \  y_3=Ta_3+a_1, \ (v_1,v_2,p-1,n_2,n_3)=(1,4,4,1,3).$$
\item $s_1<r_1$. We derive  from Theorem \ref{lem20} that $r_1-r_3+s_2-s_1<0$.  Since $vs(y_3)=\{v_1<p-1\}$  and $vs(y_2)=\{v_1<v_2\}$,  we get $$0< r_1-s_1 \le v_1-s_1\le v_1< p-1-r_3+r_1=n_3+r_1 < p-1+s_1-s_2< p-1\le 4.$$
Then $$ p=5,v_1=r_1=n_3=s_2=1, \ s_1=0, \ r_3=3 \mbox{ and } \ (v_2,n_2)=(3,1), \ \mbox{or} \  (v_2,n_2)=(4,0).$$
For   $(v_2,n_2)=(4,0)$, the obtained triple is similar to the one in $(6)$.  So, only  $(v_2,n_2)=(3,1)$  provides a new indecomposable triple. 
$$(7) \  \ \ y_2=a_1+Ta_2, \  y_3=T^3a_3+Ta_1, \ (v_1,v_2,p-1,n_2,n_3)=(1,3,4,2,1).$$ 
\ei 
$(d)$   Since $vs(y_2)=\{v_1<v_2\}$ and  $vs(y_3)=\{v_2<p-1\}$, we have  
$$ s_1\le v_1, s_2\le v_2, v_1-s_1<v_2-s_2, r_2<r_3 \mbox{ and } v_2-r_2<p-1-r_3<p-1.$$
 Then $$s_1\le v_1<v_2-s_2+s_1< v_2<p-1-r_3+r_2<p-1\le4.$$ It follows that $p=5$, $s_1=v_1=0$, $v_2=2$, $s_2=1$, $r_2=0$ and $r_3=1$. Hence we get $$(8) \ \ y_2=a_1+Ta_2, \ y_3=Ta_3+a_2, \  (v_1,v_2,p-1,n_2,n_3)=(0,2,4,1,3).$$ 
We show below that  if $p\le 5$ and $(V,T,U)$ is indecomposable, then necessarily $n_U\le 2$.\begin{lem}\label{lem22}
Let $(V,T,U)$ a triple such that $[V]=[U]$, and $\{y_1,y_2,y_3\}\inclu U$ be such that $U=vect_T\{y_1,y_2,y_3\}$. The triple $(V,T,U)$ is decomposable if one of the following condition holds,
\be
\item $\{y_1,y_2\}\inclu [V]$;
\item $y_1\in [U]$, and $(V,T,vect_T\{y_2,y_3\})$ is decomposable;
\item $Ker(T^2)\inclu U$, $T^2(y_1)=0$ and  $(V,T,vect_T\{y_2,y_3\})$ is decomposable.
\ee
\end{lem}
 \it{Proof.}
\be
\item Let $V_1=red_{y_3}$ and $V_2=W(y_3)$. We have $V=V_1\oplus V_2$, $y_3\in V_1$ and $[V]=[U]=[V_1]\oplus [V_2]$. So if
    $\{y_1,y_2\}\inclu [V]$ then $(V,T,U)$ is decomposable.
\item Suppose now $y_1\in [U]$  and that $(V,T,vect_T\{y_2,y_3\})$ is decomposable. Let $V_1$ and $V_2$ be in  $Lat(T)$ such that $V=V_1\oplus V_2$
and $vect_T\{y_2,y_3\}=(V_1\cap vect_T\{y_2,y_3\})\oplus (V_2\cap vect_T\{y_2,y_3\})$.   Since $y_1\in [U]$ and $[U]=[V_1]\oplus [V_2]$, we have 
   $(V,T,U)$ is decomposable.
\item
Also, let $V_1$ and $V_2$ in $Lat(T)$  be such that $V=V_1\oplus V_2$
and $vect_T\{y_2,y_3\}=(V_1\cap vect_T\{y_2,y_3\})\oplus (V_2\cap vect_T\{y_2,y_3\})$. We write 
$y_1=z_1+z_2$, where $z_1\in V$ and $z_2\in V_2$.  From $T^2(y_1)=0$, we obtain   $T^2(z_1)=0$  and  $T^2(z_2)=0$.    Hence $z_1\in U$ and $z_2\in U$. Thus $(V,T,U)$ is decomposable.
\ee

\begin{thm}\label{lem23} Assme that $p\le 5$.
If $(V,T,U)$ is indecomposable, then $n_U\le 2$.
\end{thm}
 \it{Proof.}

Arguing by contradiction, suppose  that $\{x_1,x_2,x_3\}$ is a basis in $[U]$ with $VP_{V,U}$. Since $(V,T,U)$ is indecomposable, we get  $n_i<v_i$ for $1\le i \le3$, $\{p-1\}\inclust vs(y_3)$ and hence $0\le n_1< v_1<p-1-1\le 3$. We distinghish  $n_2=0$ and  $n_2\neq 0$.\\
\begin{itemize}
\item $n_2=0$. We obtain $y_2=x_2$, $n_1=1$, $T(y_1)=x_1$ and  $v(y_1)=v_1-1$. Let us consider $U_1=vect_T\{y_1,y_3\}$. If $ (V,T,U_1)$ is indecomposable, then since $x_1\in U_1$, by the proof  of  Theorem \ref{lem21}, it will  follow that  $(v_1,v_2, v_3,n_1,n_3)=(1,2,4,0,2)$.  Since  $n_1=1$, we derive that $ (V,T,U_1)$ is decomposable.  By Lemma \ref{lem22} we obtain  $(V,T,U)$ is also decomposable.\\
\item $n_2\ne 0$.    Since $[V]=[U]$ then there exist $z_2\in U$ such that $T(z_2)=x_2$, and $v(z_2)=v_2-1$. We consider $ U_1=vect_T\{y_2,y_3\}$.  Arguing by contradiction again, suppose that  $ (V,T,U_1)$ is indecomposable . If $n_2=1$, then since $T(z_2)=x_2$  and $(z_2)=v_2-1$, we get $ vs(z_2)=\{v_2\}$. Hence by the proof of  Theorem \ref{lem21}, we  have $(v_1,v_2,v_3,n_2,n_3)=(1,2,4,0,3)$.  In particular $n_2=0$, which contradicts our assumption. Thus  $ n_2\ge 2$.
\end{itemize}
 In this situation, we use the proof of Theorem \ref{lem21} to get $(v_1,v_2,v_3,n_2,n_3)=(1,3,4,2,1)$. In particular we deduce that $n_3=1$.

 Now, since $T(z_3)=x_3$, and $v(z_3)=3$, we have  $vs(z_3)=\{4\}$.     It follows that $(V,T,U_1)$ is decomposable.

On the other hand,  since $card(vs(y_3))\geq 2$,    we must have $h_1\le 2$.

It is clear that if \begin{itemize}
\item $n_1=h_1$, then $(V,T,U)$ is decomposable.
\item If $n_1=0$, then by Lemma  \ref{lem22}, $(V,T,U)$ is decomposable. \item  If $n_1=1$, then from $v(z_3)=v_3-1$, $v(z_2)=v_2-1$, we get by lemma \ref{lem22},  $(V,T,U)$ is decomposable.
\end{itemize}
\begin{flushright}
$ \blacksquare$
 \end{flushright} 
From the previous sections; we have the following distribution of the $50$ indecomposable triples associated with $p\le 5$.

\begin{itemize}
\item $dim([U])=0$: Indeomposable triples are $(K^n, J_n,0)$, $1\le n\le 5$.  $N_0=5;$
\item $dim([U])= 1$: $N_1= C_1^1+C_2^1+C_3^1+C_3^3+C_4^1+C_4^3+C_5^1+C_5^3+C_5^5 =31$;
\item  $dim([U])= 2, dim([V])= 2$: $N_2= C_4^4+C_5^4=6$;
\item $dim([U])= 2, dim([V])= 3$. $N_3= 8.$ 
\end{itemize}

\end{document}